\documentclass [12pt,a4paper,reqno]{amsart}
\usepackage{amsmath,amsthm}
\usepackage{amsfonts}

\catcode`\@=11

 \def\AMSTeXfeatures{\Plainheads 
   \let\current@vert=\AMS@vert}

 \def\Plainheads{\sh@ftdiam=0.05em
   \getlabeldims
   \let\vshaftfill=\plnvsolidfill
   \let\hshaftfill=\plnhsolidfill
   \let\th@rhead=\plnrhead
   \let\th@lhead=\plnlhead
   \let\th@dnhead=\plndnhead
   \let\th@uphead=\plnuphead}
 
 \def\glet{\global\let}

 \def\LaTeXfeatures{\catcode`\@=11
   \ifx\@clnwd\undefined \nol@g
      \input ltxcode.tex \dol@g \fi
   \ltxheads \let\current@vert=\new@vert
   \providelto \catcode`\@=\active}

 \def\nol@g{\def\wlog{\edef\garbage}}
 \def\dol@g{\let\wlog=\wl@g} \let\wl@g=\wlog
 \nol@g 

 \newbox\ltobox
 \def\providelto{{\setbox\z@=
   \hbox{$\to$}\minharrlen=\wd\z@
   \global\setbox\ltobox=\hbox{$\activeat>>>$}}
   \def\lto{\mathrel{\copy\ltobox}}}

 \def\ltxheads{\sh@ftdiam=\@wholewidth
   \getlabeldims
   \let\vshaftfill= \ltxvsolidfill
   \let\hshaftfill=\ltxhsolidfill
   \let\th@rhead=\ltxrhead
   \let\th@lhead=\ltxlhead
   \let\th@dnhead=\ltxdnhead
   \let\th@uphead=\ltxuphead}
 {\catcode`\@=\active
   \gdef@#1{\csname #1\string@at\endcsname}
   \glet\activeat=@}
 \def\def@#1{\expandafter\def\csname #1@at\endcsname}

 \def@>#1>#2>{\@rrow R{#1}{#2}}
 \def@<#1<#2<{\@rrow L{#1}{#2}}
 \def@ V#1V#2V{\@rrow V{#1}{#2}}
 \def@ A#1A#2A{\@rrow A{#1}{#2}}
 \def@/#1/#2/#3/{\@rrow{#1}{#2}{#3}}
 \def@.{\ifodd\row\ifmmode\noharrow
     \else\leavevmode.\spacefactor3000 \fi
   \else\novarrow\fi}
 \def@={\ifodd\row\harrow\hequalfill{}{}%
   \else\varrow\vequalfill{}{}\fi}
 \def@:#1{\ifx=#1\harrow\deffill{}{}%
   \else\leavevmode\null:#1\fi}
 \def@|{\current@vert}
  \def\AMS@vert{\varrow\vequalfill{}{}}
  \def\new@vert#1|#2|{\ifodd\row
   \let\nextarrow\vertexvarrow
   \else\let\nextarrow\varrow\fi
   \nextarrow\vshaftfill{#1}{#2}}
 \def@-{\ifmmode\let\next\hl@ne
   \else\let\next\AMSatdash \fi \next}
  \def\hl@ne#1-#2-{\harrow\hshaftfill{#1}{#2}}
  \def\AMSatdash{\let\next\relax\leavevmode
    \def\next@{\ifx\next-%
      \def\next-{\futurelet\next\nextii@}%
     \else\def\next{\hbox{-}}\fi\next}%
    \def\nextii@{\ifx\next-\def\next-{\hbox{---}}%
      \else\def\next{\hbox{--}}\fi\next}%
    \futurelet\next\next@}
 \def@(#1){\tweenarrows{#1}}
 \def@[#1]{\setsp@n#1\relax\activeat}
 \def\fiberbox{\hbox{$\vcenter{\hr@le\hbox{\vr@le
   \kern1ex\vbox{\kern1.2ex}\vr@le}\hr@le}$}}
  \def\hr@le{\hrule height \sh@ftdiam}
  \def\vr@le{\vrule width \sh@ftdiam}
 \def@+#1+#2+#3+{\ifodd\row \harrow{#1}{#2}{#3}%
   \else \varrow{#1}{#2}{#3}\fi}

 \def\Dnarrfill{\vequalfill\Dnhe@d}
 \def\Uparrfill{\Uphe@d\vequalfill}
 
 \def\ontofill{\rtarrfill\kern-0.3em 
   \th@rhead\kern 0.3em} 

 \def\rtarrfill{\hshaftfill\th@rhead}
 \def\ltarrfill{\th@lhead\hshaftfill}
 \def\dnarrfill{\vshaftfill\th@dnhead}
 \def\uparrfill{\th@uphead\vshaftfill}
 \def\hequalfill{\plnhfill=}
 \def\deffill{:\plnhfill=}
 \def\plnvextfill#1{\setbox\z@
   \hbox{\the\textfont3 #1}%
   \dimen@=\dp\z@\advance\dimen@\ht\z@
   \copy\z@ \kern-\dimen@ 
   \cleaders\copy\z@ \vfill
   \kern-\dimen@ 
   \box\z@}
 \def\plnhfill#1{$\m@th\mkern-1.5mu\mathord#1\mkern-6mu
    \cleaders\hbox{$\mkern-2mu\mathord#1\mkern-2mu$}\hfill
    \mkern-6mu\mathord#1\mkern-1.5mu$}
 \def\vequalfill{\plnvextfill{\char'167}}
 \def\plnvsolidfill{\plnvextfill{\char'077}}
 \def\plnhsolidfill{\plnhfill-}
 \def\ltxhsolidfill{\leaders\hrule height\topofshaft depth\botofshaft
   \hfill}
 \def\ltxvsolidfill{\leaders\vrule width\sh@ftdiam\vfill}
 \def\hdashfill{\hd@sh\wd@sh
   \xleaders \hbox{\wd@sh\hd@sh\wd@sh}\hfill
   \wd@sh\hd@sh}
 \def\vdashfill{\vd@sh\wd@sh
   \xleaders \vbox{\wd@sh\vd@sh\wd@sh}\vfill
   \wd@sh\vd@sh}
 \def\dashed{\ifinmeasureCD\else
    \ifodd\row\option{\let\hshaftfill=\hdashfill}%
   \else\option{\let\vshaftfill=\vdashfill}\fi\fi}

 \newdimen\CDstrutht  \newdimen\CDstrutdp
 \newdimen\CDstrutlen \CDstrutlen=\CDstrutht
   \advance\CDstrutlen by \CDstrutdp

 \def\CDstrut{\vrule
   height \ifnum\row=1 \z@\else\CDstrutht \fi
   depth \ifnum\row=\numrows \z@ \else\CDstrutdp \fi
   width\z@}

 \newdimen\CDarrsurr \CDarrsurr=0.375em
 \newdimen\CDdashlen
 \newdimen\CDvarrlen \CDvarrlen=1.5\baselineskip
 \newdimen\minharrlen 
  \setbox\z@\hbox{$\longrightarrow$} \minharrlen=\wd\z@
 \newdimen\minCDharrlen \minCDharrlen=2.5em 
\newdimen \minc@lwd
\def\findminc@lwd{\minc@lwd=2\CDarrsurr
  \advance\minc@lwd\minCDharrlen}

 \newdimen\sh@ftdiam

 \newdimen\labelsurr \labelsurr=1.25 em

\newcount\sp@ncnt \sp@ncnt=\@ne
\newcount\sp@ncnt@ \sp@ncnt@=\@ne
\newdimen\@rrwd \newdimen\@rrdp

 \def\adjustbot#1{\option{\advance\@rrdp#1\relax}}
 
\def\pushvertex#1{\global\p@shlen#1\relax
   \global\let\maybepush=\dopush}

 \newdimen\p@shlen \p@shlen=\z@

 \let\maybepush=\relax
 \def\dopush{\ifinmeasureCD 
   \advance\locdimen by -\p@shlen 
   \else\advance \@rrwd by -\p@shlen \fi 
   \global\let\maybepush=\relax \global\p@shlen=\z@\relax}

 \def\span@ne{\global\sp@ncnt=\@ne\relax}
 \def\setsp@n#1#2{\global\sp@ncnt=#1\relax
   \ifx\relax#2\relax\else\global\sp@ncnt@=#2\relax\fi}

 \def\plnrhead{\llap{$\rightarrow\mkern-1.5mu$}}
 \def\plnlhead{\rlap{$\mkern-1.5mu\leftarrow$}}

 \def\clap#1{\hbox to \z@{\hss #1\hss}}

 \def\plndnhead{\hbox{\the\textfont3 \char'171}}
 \def\plnuphead{\hbox{\the\textfont3 \char'170}}
 \def\Dnhe@d{\hbox{\the\textfont3 \char'177}}
 \def\Uphe@d{\hbox{\the\textfont3 \char'176}}

 \def\ltxrhead{\raise\@xisheight
   \llap{\smash{\@linefnt\@getrarrow(1,0)}}}
 \def\ltxlhead{\raise\@xisheight
   \rlap{\@linefnt\@getlarrow(-1,0)}}
 \def\ltxuphead{\setbox\z@=\rlap{%
   \kern\@halfwidth\@linefnt\char'66}%
   \copy\z@\kern-\ht\z@}
 \def\ltxdnhead{\setbox\z@=\rlap{%
   \kern\@halfwidth\@linefnt\char'77}%
   \ht\z@=\z@\box\z@}

 \def\wd@sh{\kern0.5\CDdashlen}
 \def\hd@sh{\vrule height\topofshaft depth\botofshaft
    width\CDdashlen}
 \def\vd@sh{\hrule height\CDdashlen
   depth\z@ width\sh@ftdiam}

\def\xylist{14{3434}13{2414}12{1723}%
  23{1413}34{1153}11{0867}43{0707}%
  32{0580}21{0414}31{0291}41{0}}
\newcount\tgtcnt@
\def\find@xyargs{\dimen@=\@rrdp
  \advance\dimen@ by \CDstrutlen
  \tgtcnt@=\dimen@ \dimen@=\@rrwd 
  \divide\dimen@ by \@m 
  \divide \tgtcnt@ by \dimen@ 
  \expandafter\testxy\xylist\relax
  \unitlength=\@xarg\@rrdp
  \divide\unitlength by\@yarg\relax}
\def\testxy#1#2#3{\ifnum\tgtcnt@>#3
    \@xarg=#1\relax \@yarg=#2\relax
    \let\next=\ignorerest
  \else\let\next\testxy\fi\next}
\def\ignorerest#1\relax{\relax}

\let\scalefactor=\@ne
\def\SWarrow{\find@xyargs\vector
  (-\@xarg,-\@yarg)\scalefactor\hskip-\wd\@linechar}
\def\NWarrow{\find@xyargs\vector
  (-\@xarg,\@yarg)\scalefactor\hskip-\wd\@linechar}
\def\NEarrow{\find@xyargs\vector
  (\@xarg,\@yarg)\scalefactor}
\def\SEarrow{\find@xyargs\vector
  (\@xarg,-\@yarg)\scalefactor}
\def\rightupline{\find@xyargs\@linelen=\scalefactor
     \unitlength\@sline}
\def\rightdownline{\find@xyargs\@yarg=-\@yarg\relax
     \@linelen=\scalefactor\unitlength\@sline}

\def\Sim{\ifodd\row\setbox\z@=\hbox{$\sim$}\dimen@=\ht\z@
 \advance\dimen@ by -\@xisheight
  \vbox{\box\z@\kern-\@xisheight\kern\dimen@}%
  \else\hbox{$\wr$}\fi}

\def\harrow#1#2#3{\inmeasureCDtrue\findminarrwd
  {#2}{#3}{\sp@ncnt\minharrlen}\inmeasureCDfalse\span@ne
  \mathrel{\hbox{\options\hplace{#1}\ulabel{#2}\dlabel{#3}}}}

\def\noharrow{\harrow\hfill{}{}}
\def\vertexvarrow#1#2#3{\findarrdp \@rrwd=\z@ \setsp@n\@ne\@ne
  \vbox to \z@{\kern-1.2\CDstrutht
  \rlap{\options\vplace{#1}\llabel{#2}\rlabel{#3}}\vss}}

\newif\ifinmeasureCD
\def\measurelabel#1{\setbox\z@
  \hbox{$\scriptstyle#1\kern\labelsurr$}%
  \ifdim\wd\z@>\@rrwd \@rrwd=\wd\z@\fi}
\def\findminarrwd#1#2#3{\@rrwd=#3\relax
   \measurelabel{#1}\measurelabel{#2}}
\def\findCDarrwd#1#2{\@rrwd=\minCDharrlen
   \measurelabel{#1}\measurelabel{#2}%
  }

\newcount\row \row=\@ne \newcount\col \col=\@ne 
 \newcount\numrows 
\numrows=\@ne
 \newcount\numcols
\newcount\arrspan \newdimen\vrtxhalfwd  \newbox\tempbox

\def\DANABUG{\advance\col by \@ne
 \@rrwd=\minCDharrlen
  \advance\@rrwd by \vrtxhalfwd
  \advance\@rrwd by \CDarrsurr
  \ifnum\col>\numcols \numcols=\col
     \newlocdimen{col\the\col}\locdimen=\@rrwd 
  \else \ifdim\@rrwd>\c@l \c@l=\@rrwd\fi\fi}

\def\drop#1\\{
  \findvrtxhalfsum\DANABUG\advance\row by 2 \measureinit}

\def\measureinit{\col=\@ne \vrtxhalfwd=-\CDarrsurr\arrspan=\@ne\@rrwd=\z@
   \setbox\tempbox=\hbox\bgroup$}
\def\measure{
  \let\harrow\measureCDarrow
  \let\CDCR=\measureCR 
   \findminc@lwd 
  \inmeasureCDtrue
  \row=\@ne \numcols=\z@ \measureinit}

\def\endmeasure{\findvrtxhalfsum\DANABUG
  \numrows=\row 
  \inmeasureCDfalse}

\def\newlocdimen#1{\advance\dimenc@unt by \@ne
  \ifnum\dimenc@unt<\insc@unt
     \else\errmessage{No room for the CD}\fi
  \dimendef\locdimen=\dimenc@unt
  \expandafter\dimendef\csname#1\endcsname=\dimenc@unt}

 \def\r@wc@l{\csname row\the\row col\the\col\endcsname}
 \def\c@l{\csname col\the\col\endcsname}

 \def\findvrtxhalfsum{$\egroup
  \newlocdimen{row\the\row col\the\col}
  \locdimen=\vrtxhalfwd 
  \vrtxhalfwd=0.5\wd\tempbox 
  \advance\vrtxhalfwd by \CDarrsurr
  \advance\locdimen by \vrtxhalfwd 
  \advance\@rrwd by \locdimen 
  \maybepush
  \divide\@rrwd by \arrspan\relax
  \ifdim\@rrwd<\minc@lwd
    \ifnum\col>\@ne \@rrwd=\minc@lwd\fi \fi
  \loop 
    \ifnum\col>\numcols \numcols=\col
       \newlocdimen{col\the\col}
       \locdimen=\@rrwd 
    \else \ifdim\@rrwd>\c@l \c@l=\@rrwd\fi \fi
   \ifnum\arrspan>\@ne
      \advance\arrspan by -1 \advance\col by \@ne
  \repeat }

 \def\measureCDarrow#1#2#3{\findvrtxhalfsum
   \arrspan=\sp@ncnt\relax\global\sp@ncnt=1\relax
   \advance\col by \@ne
   \findCDarrwd{#2}{#3}%
   \setbox\tempbox=\hbox\bgroup$}

 \newcount\dr@tn \dr@tn=\z@
 \def\locate#1:#2{\ifinmeasureCD\else
   \count@=-#1
   \multiply\count@ by 2
   \advance\count@ by #2
   \dimen@=\count@\@rrwd
   \ifnum\dr@tn=\@ne\relax \else\dimen@=-\dimen@ \fi
   \dimen@i=\@rrdp
   \ifnum\dr@tn>\z@\advance\dimen@i by \CDstrutlen \fi
   \dimen@i=\count@\dimen@i
   \count@=#2 \multiply\count@ by 2
   \divide\dimen@ by \count@
   \divide\dimen@i by \count@
   \lift\dimen@i\nudge\dimen@\fi}

\def\betweenCDrows{\advance\row by \@ne \col=\@ne
\options}

\def\hbegin{\hbox\bgroup\kern\c@l \kern-\r@wc@l$}
\def\hend{$\glet\maybepush\relax \CDstrut\egroup}
\def\vbegin{\setbox\tempbox=\hbox\bgroup$}
\def\vend{$\egroup\ht\tempbox=\z@\dp\tempbox\CDvarrlen
  \box\tempbox}
\def\setCD{\let\harrow=\setCDarrow
  \let\CDCR=\setCR 
  \row=\@ne \col=\@ne \hbegin}
\let\endsetCD=\hend 

\def\findarrwd{\@rrwd=\z@ \count@=\col \advance\count@ by\sp@ncnt
  \loop\ifnum\count@>\col \advance\count@ by -1
      \advance\@rrwd by\csname col\the\count@\endcsname\repeat}
\def\setCDarrow#1#2#3{\kern\CDarrsurr\advance\col by \@ne
  \findarrwd \advance\@rrwd by -\r@wc@l  
  \@rrdp=\z@ 
  \maybepush
  \advance\col by -\@ne \advance\col by \sp@ncnt \span@ne
  \hbox to \@rrwd{\options
   \@rrwd=\scalefactor\@rrwd\hss
   \hplace{#1}\ulabel{#2}\dlabel{#3}\hss}%
   \kern\CDarrsurr}

\newdimen\labspacei 
\newdimen\labspaceii 

\newdimen\@xisheight
  \@xisheight=\the\fontdimen22\textfont2
\newdimen\labelskip
  \labelskip=\the\fontdimen10\textfont3 
\newdimen\topofshaft
\newdimen\botofshaft
\newdimen\botofulabel
\newdimen\topofdlabel
\def\getlabeldims{
  \topofshaft=0.5\sh@ftdiam
  \botofshaft=\topofshaft
  \advance\topofshaft by \@xisheight  
  \advance\botofshaft by -\@xisheight  
  \botofulabel=\topofshaft
  \advance\botofulabel by \labelskip
  \topofdlabel=\botofshaft
  \advance\topofdlabel by \labelskip}

\def\ulabel{\ifnum\row=\@ne\let\next\ulabeli
   \else\let\next\ulabellap\fi\next}
\def\ulabeli#1{\vbox{
  \clap{\kern-\@rrwd$\scriptstyle#1$}%
  \kern\botofulabel}\maybeoffset}
\def\ulabellap#1{\vbox to \z@{\vss
  \clap{\kern-\@rrwd$\scriptstyle#1$}%
  \kern\botofulabel}\maybeoffset}
\def\dlabel{\ifnum\row=\numrows\let\next\dlabeli
   \else\let\next\dlabellap\fi\next}
\def\dlabeli#1{\vtop{\kern\topofdlabel
  \clap{\kern-\@rrwd$\scriptstyle#1$}%
  }\maybeoffset}
\def\dlabellap#1{\vbox to \z@{\kern\topofdlabel
  \clap{\kern-\@rrwd$\scriptstyle#1$}%
  \vss}\maybeoffset}
\def\rlabel#1{\vbox to \z@{\vss
  \rlap{\kern\labelskip$\scriptstyle#1$}%
  \vss\kern-\@rrdp}\maybeoffset}
\def\llabel#1{\vbox to \z@{\vss
  \llap{$\scriptstyle#1$\kern\labelskip}%
  \vss\kern-\@rrdp}\maybeoffset}
\def\swlabel#1{\vtop{\kern0.5\@rrdp
  \llap{$\scriptstyle#1$\kern\labelskip\kern-0.5\@rrwd}
  }\maybeoffset}
\def\nwlabel#1{\vbox{
  \llap{$\scriptstyle#1$\kern\labelskip\kern-0.5\@rrwd}%
  \kern-0.5\@rrdp}\maybeoffset}
\def\selabel#1{\vtop{\kern0.5\@rrdp
  \rlap{\kern0.5\@rrwd\kern\labelskip$\scriptstyle#1$}%
  }\maybeoffset}
\def\nelabel#1{\vbox{
  \rlap{\kern0.5\@rrwd\kern\labelskip$\scriptstyle#1$}%
  \kern-0.5\@rrdp}\maybeoffset}
\def\cplace#1{\vbox to \z@{\vss
  \clap{$#1$\kern-\@rrwd}%
  \kern-\@rrdp\vss}\maybeoffset}
\def\hplace#1{\hbox to \@rrwd{#1}\maybeoffset}
\def\vplace#1{\clap{\vbox to \z@{#1\kern-\@rrdp}}\maybeoffset}

\newdimen\nudgeamount \nudgeamount=\z@
\newdimen\liftamount \liftamount=\z@
\let\maybeoffset\relax
\newbox\offsetbox \newdimen\lastheight
\def\dooffset{
  \setbox\offsetbox=\lastbox \lastheight=\ht\offsetbox 
  \setbox\offsetbox=\vbox{\kern-\liftamount\box\offsetbox}%
  \ht\offsetbox=\lastheight
  \kern\nudgeamount\box\offsetbox\kern-\nudgeamount
  \global\nudgeamount=\z@ \global\liftamount=\z@
  \glet\maybeoffset=\relax}
\def\nudge#1{\ifinmeasureCD\else
  \global\advance\nudgeamount#1\relax
  \global\let\maybeoffset\dooffset\fi}
\def\lift#1{\ifinmeasureCD\else
  \global\advance\liftamount#1\relax
  \global\let\maybeoffset\dooffset\fi}

\def\findarrdp{\@rrdp=\CDvarrlen
  \ifnum\sp@ncnt@>1
    \advance\@rrdp by \CDstrutlen
    \multiply\@rrdp by \sp@ncnt@
    \advance\@rrdp by -\CDstrutlen \fi
 }

\def\varrow#1#2#3{\ifnum\sp@ncnt>\@ne 
     \sp@ncnt@=\sp@ncnt\relax\fi
  \findarrdp \@rrwd=\z@ 
  \kern\c@l
   \hbox to \z@{\options
   \@rrdp=\scalefactor\@rrdp
    \hss\vplace{#1}\llabel{#2}\rlabel{#3}\hss}%
  \global\advance\col by \@ne \setsp@n\@ne\@ne
  }

\def\novarrow{\varrow\vfill{}{}}

\def\tweenarrows#1{\findarrwd \findarrdp \setsp@n\@ne\@ne
  \rlap{\options\cplace{#1}}}

\def\usarrow #1#2#3{\dr@tn=\@ne
  \findarrwd \findarrdp \setsp@n\@ne\@ne 
  \rlap{\options\cplace{#1}\nwlabel{#2}\selabel{#3}}%
  \dr@tn=\z@}
\def\dsarrow #1#2#3{\dr@tn=\tw@
  \findarrwd \findarrdp \setsp@n\@ne\@ne 
  \rlap{\options\cplace{#1}\swlabel{#2}\nelabel{#3}}%
  \dr@tn=\z@}
 \def\@rrow#1{\csname #1@rrow\endcsname}
 \def\R@rrow{\harrow \rtarrfill}
 \def\L@rrow{\harrow \ltarrfill}
 \def\V@rrow{\varrow \dnarrfill}
 \def\A@rrow{\varrow \uparrfill}
 \def\SE@rrow{\dsarrow \SEarrow}
 \def\NW@rrow{\dsarrow \NWarrow}
 \def\SW@rrow{\usarrow \SWarrow}
 \def\NE@rrow{\usarrow \NEarrow}
 \def\DS@rrow{\dsarrow \dnslope}
 \def\US@rrow{\usarrow \upslope}
 \def\upslope{\find@xyargs
       \@linelen=\unitlength\@sline}
 \def\dnslope{\find@xyargs\@yarg=-\@yarg\relax
       \@linelen=\unitlength\@sline}

\newtoks\optionlist 
\optionlist={}
\let\options\relax
\def\dooptions{\the\optionlist\global\optionlist={}%
  \glet\options=\relax}
\def\option#1{\ifinmeasureCD\else
  \glet\options=\dooptions
  \global\optionlist=\expandafter{\the\optionlist\relax#1}\fi}
\def\wider#1{\ifinmeasureCD\else
   \option{\advance\@rrwd by #1}\fi}
\def\deeper#1{\ifinmeasureCD\else
   \option{\advance\@rrdp by #1}\fi}

{\def\\{\global\let\sptoken= }\\ }

\def\CR{\futurelet\nexttok\testCR}
\def\testCR{\ifx\nexttok\sptoken
   \let\next\eatspaceCR\else\let\next\CDCR\fi\next}
\def\eatspaceCR#1 {\CR}
\def\measureCR{\ifx\nexttok\endmeasure\let\nextCR\relax
    \else\let\nextCR\drop\fi\nextCR}
\def\setCR{\ifodd\row
  \ifx\nexttok\endsetCD\else\hend\betweenCDrows\vbegin\fi
  \else\vend\betweenCDrows\hbegin\fi}

\countdef\dimenc@unt=11
\def\CD#1\endCD{
   \begingroup\let\\=\CR
  \m@th\offinterlineskip
   \measure#1\endmeasure\null\,\vcenter{\setCD#1\endsetCD}\,
   \endgroup
    }

\ifx\@clnwd\undefined \nol@g\else\catcode`\ =14\relax\fi
 \font\@linefnt=line10 
 \newcount\@tempcnta
 \newcount\@tempcntb
 \newdimen\@tempdima
 \newdimen\@tempdimb
 \newdimen\@wholewidth
 \newdimen\@halfwidth
   \@wholewidth\fontdimen8\@linefnt \@halfwidth .5\@wholewidth
 \newdimen\unitlength
 \newcount\@xarg
 \newcount\@yarg
 \newcount\@yyarg
 \newbox\@linechar
 \newdimen\@linelen
 \newdimen\@clnwd
 \newdimen\@clnht
 \newif\if@negarg
 
 \def\@whilenoop#1{}

 \def\@whiledim#1\do #2{\ifdim #1\relax#2\@iwhiledim{#1\relax#2}\fi}

 \def\@iwhiledim#1{\ifdim #1\let\@nextwhile=\@iwhiledim 
         \else\let\@nextwhile=\@whilenoop\fi\@nextwhile{#1}}

 \def\@sline{\ifnum\@xarg< 0 \@negargtrue \@xarg -\@xarg \@yyarg -\@yarg
   \else \@negargfalse \@yyarg \@yarg \fi
 \ifnum \@yyarg >0 \@tempcnta\@yyarg \else \@tempcnta -\@yyarg \fi
 \ifnum\@tempcnta>6 \@badlinearg\@tempcnta0 \fi
 \ifnum\@xarg>6 \@badlinearg\@xarg 1 \fi
 \setbox\@linechar\hbox{\@linefnt\@getlinechar(\@xarg,\@yyarg)}%
 \ifnum \@yarg >0 \let\@upordown\raise \@clnht\z@
    \else\let\@upordown\lower \@clnht \ht\@linechar\fi
 \@clnwd=\wd\@linechar
 \if@negarg \hskip -\wd\@linechar \def\@tempa{\hskip -2\wd\@linechar}\else
      \let\@tempa\relax \fi
 \@whiledim \@clnwd <\@linelen \do
   {\@upordown\@clnht\copy\@linechar
    \@tempa
    \advance\@clnht \ht\@linechar
    \advance\@clnwd \wd\@linechar}%
 \advance\@clnht -\ht\@linechar
 \advance\@clnwd -\wd\@linechar
 \@tempdima\@linelen\advance\@tempdima -\@clnwd
 \@tempdimb\@tempdima\advance\@tempdimb -\wd\@linechar
 \if@negarg \hskip -\@tempdimb \else \hskip \@tempdimb \fi
 \multiply\@tempdima \@m
 \@tempcnta \@tempdima \@tempdima \wd\@linechar \divide\@tempcnta \@tempdima
 \@tempdima \ht\@linechar \multiply\@tempdima \@tempcnta
 \divide\@tempdima \@m
 \advance\@clnht \@tempdima
 \ifdim \@linelen <\wd\@linechar
    \hskip \wd\@linechar
   \else\@upordown\@clnht\copy\@linechar\fi}
 
 \def\@getlinechar(#1,#2){\@tempcnta#1\relax\multiply\@tempcnta 8
 \advance\@tempcnta -9 \ifnum #2>0 \advance\@tempcnta #2\relax\else
 \advance\@tempcnta -#2\relax\advance\@tempcnta 64 \fi
 \char\@tempcnta}
 
 \def\vector(#1,#2)#3{\@xarg #1\relax \@yarg #2\relax
 \@tempcnta \ifnum\@xarg<0 -\@xarg\else\@xarg\fi
 \ifnum\@tempcnta<5\relax
 \@linelen=#3\unitlength
 \ifnum\@xarg =0 \@vvector 
   \else \ifnum\@yarg =0 \@hvector \else \@svector\fi
 \fi
 \else\@badlinearg\fi}
 
 \def\@svector{\@sline
 \@tempcnta\@yarg \ifnum\@tempcnta <0 \@tempcnta=-\@tempcnta\fi
 \ifnum\@tempcnta <5
   \hskip -\wd\@linechar
   \@upordown\@clnht \hbox{\@linefnt  \if@negarg 
   \@getlarrow(\@xarg,\@yyarg) \else \@getrarrow(\@xarg,\@yyarg) \fi}%
 \else\@badlinearg\fi}
 
 \def\@getlarrow(#1,#2){\ifnum #2 =\z@ \@tempcnta='33\else
 \@tempcnta=#1\relax\multiply\@tempcnta \sixt@@n \advance\@tempcnta
 -9 \@tempcntb=#2\relax\multiply\@tempcntb \tw@
 \ifnum \@tempcntb >0 \advance\@tempcnta \@tempcntb\relax
 \else\advance\@tempcnta -\@tempcntb\advance\@tempcnta 64
 \fi\fi\char\@tempcnta}
 
 \def\@getrarrow(#1,#2){\@tempcntb=#2\relax
 \ifnum\@tempcntb < 0 \@tempcntb=-\@tempcntb\relax\fi
 \ifcase \@tempcntb\relax \@tempcnta='55 \or 
 \ifnum #1<3 \@tempcnta=#1\relax\multiply\@tempcnta
 24 \advance\@tempcnta -6 \else \ifnum #1=3 \@tempcnta=49
 \else\@tempcnta=58 \fi\fi\or 
 \ifnum #1<3 \@tempcnta=#1\relax\multiply\@tempcnta
 24 \advance\@tempcnta -3 \else \@tempcnta=51\fi\or 
 \@tempcnta=#1\relax\multiply\@tempcnta
 \sixt@@n \advance\@tempcnta -\tw@ \else
 \@tempcnta=#1\relax\multiply\@tempcnta
 \sixt@@n \advance\@tempcnta 7 \fi\ifnum #2<0 \advance\@tempcnta 64 \fi
 \char\@tempcnta}
\catcode`\ =10

\dol@g 
\catcode`\@=\active
\LaTeXfeatures

 \input amssymb.sty
\DeclareMathOperator{\Com}{Com}
\DeclareMathOperator{\Ob}{Ob}
\DeclareMathOperator{\Grp}{Grp}
\DeclareMathOperator{\Hom}{Hom}
\DeclareMathOperator{\Ass}{Ass}
\DeclareMathOperator{\QSC}{QSC}

\DeclareMathOperator{\LSC}{LSC}

\DeclareMathOperator{\End}{End}
\DeclareMathOperator{\Set}{Set}
\DeclareMathOperator{\Int }{Int}
\DeclareMathOperator{\Lie }{Lie}

\newcommand{\ef}{\end{equation}}
\chardef\bslash=`\\ 

\newtheorem{thm}{Theorem}
\newtheorem*{thm*}{Theorem}
\newtheorem*{conjecture*}{Conjecture}

\newtheorem{cor}{Corollary}

\newtheorem{prop}{Proposition}
 \newtheorem{prob}{Problem}

\theoremstyle{definition}
\newtheorem{defn}{Definition}
\newtheorem*{remark*}{Remarks}
\newtheorem*{examples*}{Examples}
\newtheorem*{defn*}{Definition}
\newtheorem*{cor*}{Corollary}

\newcommand{\thmref}[1]{Theorem~\ref{#1}}
\newcommand{\secref}[1]{Section~\ref{#1}}

\newcommand{\propref}[1]{Proposition~\ref{#1}}

\numberwithin{equation}{section}

\newcommand{\Th}{\Theta}

 \renewcommand{\sectionmark}[1]{}

\newcommand{\Cl} {\operatorname{Cl}}

\newcommand{\wedgel}{\operatornamewithlimits{\bigwedge}\limits}
\newcommand{\capl}{\operatornamewithlimits{\bigcap}\limits}

 \usepackage{amsfonts}

 \newcommand{\dl}{\delta}

\newcommand{\ov}{\overline}
\newcommand{\vp}{\varphi}

 \newcommand{\Mod}{\operatorname{Mod}}
\newcommand{\Var}{\operatorname{Var}}

\newcommand{\Aut}{\operatorname{Aut}}
\newcommand{\Ker}{\operatorname{Ker}}

 \date{}
\begin{document}

\title{Algebras with the same (algebraic) geometry}
\author[B. Plotkin]{B. Plotkin\\
 Institute of Mathematics \\
 Hebrew University, 91803 Jerusalem, Israel}


\maketitle

\begin{abstract}
Some basic notions of classical algebraic geometry can be defined
in arbitrary varieties of algebras $\Theta.$ For every algebra $H$ in
$\Theta$ one can consider algebraic geometry in $\Theta$ over $ H.$
Correspondingly, algebras in $\Theta$ are considered with the
emphasis on equations and geometry. We give examples of geometric
properties of algebras in $\Theta$ and of geometric relations
between them. The main problem considered in the paper is when
different $H_1$ and $H_2$ have the same geometry.
\end{abstract}

\baselineskip 20pt
\bigskip

\centerline{INTRODUCTION}\label{Intro}

\subsection{}
Speaking on universal algebraic geometry, we assume that the basic
variety is arbitrary or sufficiently arbitrary. Under
non-classical algebraic geometry we mean algebraic geometry in
various specific (fixed) varieties $\Theta$, i.e., non-classical
stands for not necessarily classical.
One can consider algebraic geometry in
groups, in rings (associative or Lie), and in other structures.
All this is united by the general idea of nonclassical algebraic
geometry. Hence, there appeared universal problems and problems
arising from the peculiarities of a concrete variety $\Th.$

We distinguish varieties $\Com$-$P$, $\Ass$-$P$ and $\Lie$-$P$.
The first one is the variety of all commutative and associative
algebras with the unit over the field $P$. The geometry associated
with this variety is regarded as the classical algebraic geometry
over $P$. The second one is the variety of all associative (not
necessarily commutative) algebras with the unit over $P.$
$\Lie$-$P$ is the variety of all Lie algebras over $P$.

For every algebra $H\in\Th$ we have its algebraic structure, its
logic and its geometry. The interaction of these three components
is the main idea of the theory under consideration. This leads to
a number of new problems. For example, when do two algebras $H_1$
and $H_2$ have the same geometry, and how does one understand this
fact. ``The same algebra" means isomorphism of algebras, ``the
same logic" can be treated as the coincidence of elementary
theories.

\subsection{} With every algebra $H\in\Th$ we associate two categories.
They are the category $K_\Theta(H)$ of algebraic sets over $H$,
and the category of algebraic varieties $\tilde K_\Th(H)$. Here,
algebraic variety is viewed as an algebraic set, considered up to
isomorphism of algebraic sets. Thus, the category $\tilde
K_\Th(H)$ is the skeleton of the category $K_\Th(H).$ Both
categories represent the geometry of $H$ and are geometrical
invariants of the algebra $H.$ Now we can view geometries in the
algebras $H_1$ and $H_2$ to be the same if the categories
$K_\Theta(H_1)$ and $K_\Theta(H_2)$ are isomorphic or the
categories $\tilde K(H_1)$ and $\tilde K_\Th(H_2)$ are isomorphic.
On the other hand, in the category theory it is known that two
categories have isomorphic skeletons if and only if these
categories are equivalent. Hence, we distinguish two problems:
\begin{enumerate}
\item[1)] When are the categories $K_\Th(H_1)$ and $K_\Th(H_2)$ isomorphic?
\item[2)] When are these categories equivalent?
\end{enumerate}
\noindent
 In fact, we consider the formulated problems in respect to special correct isomorphism
and correct equivalence. This approach reflects the idea of
coincidence of geometries. Correctness is inspired by the essence
of the matter, as explained later.

Let us present here two specific results.

The first one relates to classical algebraic geometry over a field
$P.$ For an arbitrary extension $L$ of the field $P$ denote by
$K_P(L)$ the corresponding category $K_\Th(L).$

Let $L_1,$ $L_2$ be two extensions of the field $P.$
The following conditions are equivalent:
\begin{enumerate}
\item[1)] The categories $K_P(L_1)$ and $K_P(L_2)$ are correctly isomorphic.
\item[2)] These categories are correctly equivalent.
\item[3)] There exists an extension $L$ of the field $P$ such that $L_1$ and $L$ are
semi-isomorphic, and $L_2$ and $L$ have the same quasi-identities.
\end{enumerate}

The second, more simple, result relates to groups.

Let $\Th=\Grp $ be the variety of all groups, $H_1$ and $H_2$ two nonperiodical abelian groups.
Then the following three conditions are equivalent:
\begin{enumerate}
\item[1)] The categories $K_\Th(H_1)$ and $K_\Th(H_2)$ are correctly isomorphic.
\item[2)] These categories are correctly equivalent.
\item[3)]   $H_1$ and $H_2$   have the same quasi-identities.
\end{enumerate}

Other cases of groups and algebras are considered in the same
spirit.

Let us note that correctness of the isomorphism of categories
$K_\Th(H_1)$ and $K_\Th(H_2)$ is well coordinated with the
lattices of algebraic sets in affine spaces.

Any category $K_\Th(H),$  can be considered from the point of view
of the possibility to solve systems of equations in the algebra
$H$ . This category is,  in some sense, a measure of algebraic
closeness of the given $H$, depending on the structure of algebra
$H$.

An important part in the proofs is played by investigation of
automorphisms of categories of free algebras of varieties. For any
$\Th$ denote by $\Th^0$ the category of free in $\Th$ algebras
$W=W(X)$ with $X$ finite. Automorphisms and autoequivalencies of
such a category $\Th^0$ are essentially tied with the geometry in
$\Th.$

\subsection{}
We consider also the category $K_\Th$ of algebraic sets over different
$H\in\Th.$
Its skeleton $\tilde K_\Th$ is a category of algebraic varieties over
different $H.$
Both these categories are geometrical invariants of the whole variety $\Th.$
There naturally arise problems on isomorphism and equivalence
of different $K_{\Th_1}$ and
$K_{\Th_2}.$
Here $\Th_1$ and $\Th_2$ could be subvarieties of some big variety $\Th.$
\subsection{}
Let us make some notes on the plan of the paper. The paper is
organized as follows. At the beginning we recall basic
definitions. The second, third and fourth sections are devoted to
special notions which play main part in the solution of the
problem whether geometries in different algebras are the same.

In the fifth section we present universal theorems on coincidence
of geometries.  In the seventh section these universal theorems
are specialized for varieties $\Com$-$P$, $\Ass$-$P$ and
$\Lie$-$P$. The previous sixth section contains preparation
material for the final seventh one.

\section{Basic definitions}\label{definitions}
\subsection{Algebraic sets and closed congruences}

General results about geometry in groups can be found in
[2,3,4,16,17,18,36,37]. Fix an arbitrary variety of algebras
$\Th.$ Denote by $\Th^0$ the category of free in $\Th$ algebras
$W=W(X)$ with $X$ fixed. All $X$ are supposed to be subsets of
some infinite universum $X^0.$ Thus, $\Th^0$ is a small category.

Fix an algebra $H$ in $\Th.$ Consider the set of homomorphisms
$\Hom(W,H)$ as an affine space over $H.$ Here points are
homomorphisms $\mu: W(X)\to H.$ For $X=\{x_1,\dots, x_n\}$, we
have a bijection $\alpha_X:\Hom(W,H)\to G^{(n)}$ by the rule
$\alpha_X(\mu)=(\mu(x_1),\dots,\mu(x_n)).$ If, further, $(w,w')$
is a pair of elements in $W,$ then the point $\mu:W\to H$ is a
solution of the equation $w=w'$ in $H$ if $w^\mu=w'{}^\mu,$\
$(w,w')\in\Ker\mu.$ The kernel of a homomorphism is a congruence
of the algebra $W.$

Let now $A$ be a subset in the affine space $\Hom(W,H)$ (a set of points), and $T$ a binary
relation in $W$ (a set of pairs $(w,w');w,w'\in W).$
We set:
\begin{equation*}
\begin{cases}
T'=T_H'=A=\{\mu:W\to H\bigm|  T\subset \Ker\mu\}\\
A'=A_W'=T=\capl_{\mu\in A}\Ker\mu
\end{cases}
\end{equation*}
This gives the Galois correspondence between sets of points and
binary relations. We call a set of points $A$ such that $A=T'$ for
some $T$ an algebraic (closed) set in $\Hom(W,H).$ A relation $T$
with $T=A'$ for some $A$ is a congruence in $W.$ We call such a
congruence an $H$-closed one.

For every $A$ we have a closure $A_H''=(A')_H'$ and
$T_H''=(T_H')_W'$ holds for every $T$. It is easy to understand
that the congruence $T$ in $W$ is $H$-closed if and only if there
is an injection $W/T\to H^I$ for some $I.$

We need some auxiliary notions.

The class of algebras $\frak X \subset\Th$ is called a {\it
prevariety} if $\frak X$ is closed under cartesian products and
subalgebras. For an arbitrary $\frak X$ the  corresponding closure
up to prevariety is $SC(\frak X).$ Here $S$ and $C$ are  closure
operators on classes of algebras: $C$ under cartesian products and
$S$ under subalgebras. We can now say that the congruence $T$ in
$W$ is $H$-closed if and only if $W/T\in SC(H).$ Besides, if $T$
is an arbitrary binary relation in $W,$ then $T_H''$ is an
intersection of all congruences $T_\alpha$ with $T\subset
T_\alpha$ and $W/T_\alpha\in SC(H).$

Consider further formulas of the form
\begin{equation}\label{$*$}
\left(\wedgel_{(w,w')\in T}(w\equiv w')\right)\Rightarrow (w_0,w_0'). \tag{$*$}
\end{equation}
We call them {\it generalized (infinitary) quasi-identities}. If
$T$ is finite, we have an ordinary quasi-identity.

The following Proposition easily follows from the definitions.

\begin{prop}\label{prop1}
The inclusion $(w_0,w_0')\in T_H''$ takes place if and only if the formula \eqref{$*$} holds in
the algebra $H.$
\end{prop}

\subsection{Categories of algebraic sets}

Define a category of affine spaces $K_\Th^0(H).$
Objects of this category have the form $\Hom(W,H) $ where $W$ is an object of the category
$\Th^0.$
Morphisms
$$\tilde s: \Hom(W(X),H)\to\Hom(W(Y),H)$$
are given by the morphisms in $\Th^0$
$$s:W(Y)\to  W(X)$$
by the rule ${\tilde s}(\nu)=\nu s$ for every $\nu: W(X)\to H.$ We
have here a contravariant functor
$$\Th^0\to K_\Th^0(H).$$
This functor determines duality of categories, if the algebra $H$
generates the whole variety, i.e., $\Th=\Var(H)=\QSC(H).$ Here,
$Q$ is the operator of taking  homomorphic images.

Let us now define the category $K_\Th(H)$ of all algebraic sets
over $H.$ Objects of this category have the form $(X,A)$, where
$A$ is an algebraic set in the affine space $\Hom(W(X),H).$ We
consider the affine space $\Hom(W(X),H)$ as an object of the
category $K_\Th(H)$ as well. As an algebraic set it is determined
by one equation $x=x.$ Morphisms $[s]: (X,A)\to (Y,B)$ are
determined by homomorphisms $s: W(Y)\to W(X)$ with the property
$\tilde s(\nu)=\nu s\in B$ if $\nu\in A.$ These are exactly those
$s$ for which $(w^s,w'{}^s)\in A'$ if $(w,w')\in B'.$

For such $s$ we have a homomorphism
$$\ov s: W(Y)/B'\to W(X)/A'.$$
Simultaneously we have a mapping $[s]: A\to B$ and consider it as a morphism in the category
$K_\Th(H).$

It is clear that the category of affine spaces $K_\Th^0(H)$ is a
subcategory of the category $K_\Th(H).$

Let us define the category $C_\Th(H).$ Its objects are of the form
$W(X)/T,$ where $W=W(X)$ is an object of the category $\Th^0$ and
$T$ is an $H$-closed congruence in $W(X).$  Morphisms in
$C_\Th(H)$ are homomorphisms of algebras, and the category
$C_\Th(H)$ is a full subcategory of the category $\Th.$

Let us note further that if $\Var(H)=\Theta$ the transition
$(X,A)\to W(X)/A'$ determines duality of categories $K_\Th(H)$ and
$C_\Th(H).$

We denote the corresponding skeletons of categories by $\tilde
K_\Th(H)$ and $\tilde C_\Th(H).$ These two categories are also
dual. Objects of the category $\tilde K_\Th(H)$ are called
algebraic varieties over the algebra $H.$ They are algebraic sets
over $H$ considered up to isomorphisms in $K_\Th(H).$

The following proposition \cite{29} takes place.

{\it Let $\Th_1=\Var(H).$
Then the categories $K_\Th(H)$ and $K_{\Th_1}(H)$ are isomorphic.
The categories $C_\Th(H)$ and $C_{\Th_1}(H)$ are isomorphic as well.}

For definitions of the categories $K_\Th$ and $C_\Th$ see \cite{33}, \cite{28}.

\section{Geometrical equivalence of algebras}\label{Geoequiv}

\subsection{Preliminaries}

\begin{defn}
Algebras $H_1$ and $H_2$ in $\Th$ are called geometrically
equivalent if $T_{H_1}''=T_{H_2}''$ in $W$ for any $W=W(X)$ and
$T$ in $W$.
\end{defn}

\propref{prop2} follows from the  \propref{prop1}.

\begin{prop}\label{prop2}
The algebras $H_1$ and $H_2$ are geometrically equivalent if and
only if their generalized quasi-identities coincide.
\end{prop}

This implies

\begin{prop}\label{prop3}
If $H_1$ and $H_2$ are geometrically equivalent  then they have
the same quasi-identities, i.e., the quasivarieties $q\Var(H_1)$
and $q\Var(H_2)$ coincide.
\end{prop}

The converse statement is not valid in general case (see Theorem
2).

For every class of algebras $\frak X\subset\Th$, we define a class $L\frak X$ as follows: $H\in
L\frak X$ if every finitely-generated subalgebra $H_0$ in $H$ belongs to the class $\frak X.$

The class $\LSC(\frak X)$ is a locally-closed prevariety, generated by the class $\frak X.$

A.I. Maltsev [21], [22] proved that if the class $SC(\frak X)$ is
axiomatizable, then this class is a quasivariety. In this case,
$LSC(\frak X)$ is a quasivariety as well. In the general case,
such a class is not axiomatizable. However, the following theorem
takes place \cite{34}:

{\it For any $\frak X$ the class $LSC(\frak X)$ is determined by
infinitary quasi-identities of the class $\frak X.$}

Here arises a natural question what are the  conditions  providing
$LSC(\frak X)=q\Var(\frak X)$. This question is related to
\propref{prop3}.

We call the class $\frak X$  {\it logically compact
$(q_\omega$-compact \cite{27})} if each of its infinitary
quasi-identity $\left(\wedgel_{(w,w')\in T}(w\equiv
w')\right)\Rightarrow w_0\equiv w_0',$ where $T$ is a binary
relation in $W=W(X)$, reduced to a finite quasi-identity
$\Bigg(\wedgel_{(w,w')\in T_0}(w\equiv w')\Bigg)\Rightarrow
w_0=w_0'$ with a finite subset $T_0$ in $T.$

We can now claim that if $\frak X$ is a logically compact class, then
$$LSC(\frak X)=q\Var(\frak X).$$
Actually, the opposite is also true (see below).

Let us  note also that the problem of coincidence of classes
$SC(\frak X)$ and $q\Var(\frak X)$ was solved by V.A. Gorbunov
\cite{13}.

\begin{prop}\label{prop4} \cite{35}
Algebras $H_1$ and $H_2$ are geometrically equivalent if and only if
$$LSC(H_1)=LSC(H_2).$$
\end{prop}

\subsection{Geometrically noetherian algebras}

We introduce, first, some new definitions.

\begin{defn}
An algebra $H\in\Th$ is called geometrically noetherian if for any
$W=W(X)$ and $T$ in $W$ there exists a finite subset $T_0$ of $T$
such that $T_H''=T_{0H}''.$
\end{defn}

The following proposition is proved in a standard way.

\begin{prop}\label{prop5}
An algebra $H\in\Th$ is geometrically noetherian if and only if
for every $W\in\Ob\Th^0$ the ascending chain condition for
$H$-closed congruences holds.
\end{prop}

The equivalent condition is descending chain condition for
algebraic sets in $\Hom(W,H)$ for every $W\in\Ob\Th^0.$

\begin{defn}
We call a variety $\Th$ noetherian if every $W\in\Ob\Th^0$ is noetherian (by congruences).

Obviously, if $\Th$ is a noetherian variety then every algebra $H\in\Th$ is geometrically
noetherian.
\end{defn}

\begin{examples*}
 \vskip .05in
\begin{enumerate}
\item[1)] A classical variety $\Com$-$P$ is noetherian.
\item[2)] All noetherian subvarieties in $\Ass$-$P$ are described \cite{1}.
\item[3)] The variety $\frak N_c$ of all nilpotent groups of the nilpotency class $c$ is
noetherian.
\item[4)] Every variety consisting of locally finite groups is noetherian.
\item[5)] A variety of the form $\frak N_c\Th,$ where $\Th$ is a locally finite variety, is
noetherian.
\item[6)] A free group $F(X)$ with finite $X$ is geometrically noetherian \cite{15}.
\item[7)] Finitely-dimensional associative and Lie algebras are geometrically noetherian
\cite{5}.
\end{enumerate}
\end{examples*}

We note that the notion of geometrical noetherianity of an
algebra, as well as the notion of geometrically equivalence of
algebras, does not depend on the choice of variety containing the
algebras under consideration.

Now let us generalize the notion of geometrical noetherianity.

\subsection{Logical noetherianity}

\begin{defn}
An algebra $H\in\Th$ is called locally geometrically noetherian if
for every free algebra $W$ and every set $T$ in $W$ and for every
pair $(w_0,w_0')\in T_H''$ there exists a finite subset $T_0$ in
$T,$ depending, generally, on $(w_0,w_0')$, such that
$(w_0,w_0')\in T_{0H}''.$
\end{defn}

\propref{prop6} follows directly from above.

\begin{prop}\label{prop6}
The algebra $H\in\Th$ is locally geometrically noetherian if every
infinitary quasi-identity in $H$ is reduced in $H$ to a finite
quasi-identity.
\end{prop}

This is the reason why locally geometrically noetherian algebras
we call also logically noetherian.

\begin{prop}\label{prop7}
The algebra $H$ is logically noetherian if and only if the union of any directed system of
$H$-closed congruences is also an $H$-closed congruence for every $W\in\Ob\Th^0.$
\end{prop}

\begin{proof}
Let the algebra $H$ be logically noetherian and $T $ a union of some directed system of
$H$-closed congruences $T_\alpha,$\ $\alpha\in I.$
$T$ is a congruence.
We need to check that it is $H$-closed.

Take $T_H''$ and let it contain the pair $(w,w').$ Find a finite
subset $T_0$ in $T$ with $(w,w')\in T_{0H}''.$ We have $T_\alpha$
with $T_0\subset T_\alpha.$ Then $(w,w')\in T_{0H}''\subset
T_{\alpha H}''=T_\alpha\subset T.$ Thus $(w,w')\in T,$\ $T=T_H''.$

To prove the opposite, assume the condition of directed systems of $H$-closed congruences.

Take an infinite set $T$ in $W.$ Consider in $T$ all  possible
finite subsets $T_\alpha.$ All $T_{\alpha H}''$ constitute a
directed system of $H$-closed congruences. Let $T_1$ be the union
of all congruences of this system. $T\subset T_1\subset T_H''.$
Since $T_1$ is $H$-closed, then $T_1=T_H''.$ If $(w,w')\in
T_H''=T_1,$ then $(w,w')\in T_{\alpha H}''$ for some $\alpha.$
This means that the algebra $H$ is logically noetherian.

It is clear that geometrical noetherianity of algebras implies its
logical noetherianity. Show that the opposite is not true for the
case   of groups $\Th=\Grp.$ Consider free groups $F=F(X),$\ where
$X$ are finite subsets in $X^0,$ and consider all possible
invariant subgroups $U$ in them. Denote by $H$ the discrete direct
product of all $F(X)/U.$ We have injections $F(X)/U\to H.$
Therefore, all invariant subgroups in every $F(X)$ are $H$-closed.
From  this it follows that the group $H$ is not geometrically
noetherian. However,  it is logically  noetherian  by
\propref{prop7}.
\end{proof}

Similar examples can be found in the variety $\Ass$-$P$ and various other cases.

\subsection{Logical notherianity and geometrical equivalence}

\begin{thm}\label{thm1} \cite{27}
The equality $LSC(H)=q\Var(H)$ takes place
if and only if the algebra $H$ is logically noetherian.
\end{thm}

\begin{proof}
In \cite{27} the theorem is proved for groups, but similar
considerations are valid in the general situation. Note that
\cite{27} uses a different term ($q_\omega$ - compactness) instead
of the term ``logical noetherianity."

We present the  proof for arbitrary $\Th,$ taking into account, in
particular, the case of associative and Lie algebras.

Note first of all that the algebra $H$ is logically noetherian if the class $\frak X,$
 consisting
of one algebra $H$, is logically compact.
Thus, if the algebra $H$ is logically noetherian, then $LSC(H)=q\Var(H).$

We now prove the opposite. Let $LSC(H)=q\Var(H)$ be given. Check
that the algebra $H$ is logically noetherian. Take an algebra
$W=W(X)\in\Ob\Th^0.$ Take a congruence $T$ in $W$ which is the
union of the directed system of $H$-closed congruences
$T_\alpha,$\ $\alpha\in I.$ We want to verify that $T$ is
$H$-closed as well, i.e., $W/T\in SC(H).$ In our conditions we
just need to check that every quasi-identity of the algebra $H$
holds in $W/T.$

Let the quasi-identity
\begin{equation}\label{*****}
w_1\equiv w_1'\wedge\dots\wedge w_n\equiv w_n'\to w_0\equiv w_0'\tag{$**$}
\end{equation}
be written in elements from $W(Y),$ and let it be fulfilled in the
algebra $H.$ Check that it holds also in $W/T.$

 Take   a
homomorphism $\mu: W(Y)\to W(X)/T,$ and the corresponding
commutative diagram
$$
\CD
W(Y) @>\mu_0>> W(X)\\
@. @/SE/ \mu  // @VV \nu V\\
@. W(X)/T \\
\endCD
$$
 Here $\nu$ is  a natural homomorphism.
Besides, for every $\alpha\in I$ consider natural homomorphisms $\nu_\alpha: W(X)\to
W(X)/T_\alpha.$
Assume that $w_i^\mu=w_i'{}^\mu;$ \ $w_i^{\mu_0\nu}=w_i'{}^{\mu_0\nu}$ holds for every
$i=1,\dots,n.$
We can choose $\alpha\in I$ such that $w_i^{\mu_0\nu_\alpha}=w_i'{}^{\mu_0\nu_\alpha}.$
We proceed from the homomorphism $\nu_\alpha\mu_0:W(Y)\to W(X)/T_\alpha.$

Since the quasi-identity \eqref{*****} holds in $W(X)/T_\alpha,$
we have also $w_0^{\mu_0\nu_\alpha}=w_0'{}^{\mu_0\nu_\alpha}.$
 The last formula gives
$ w_0^{\mu_0\nu}=w_0'{}^{\mu_0\nu};$\ $w_0^\mu=w_0'{}^\mu.$ This
means that the quasi-identity \eqref{*****}  holds in the algebra
$W(X)/T$ and the congruence $T$ is $H$-closed. Hence the algebra
$H$ is logically noetherian.
\end{proof}

Note that similar arguments can be used in the case when instead
of one algebra $H$ we take an arbitrary logically compact class.
(see \cite{34}).

\thmref{thm2} easily follows from the theorem just proved (see \cite{27}).

\begin{thm}\label{thm2}
If the algebra $H=H_1\in\Th$ is not logically noetherian, then
there exists its ultrapower $H_2$ which  is not geometrically
equivalent to the algebra $H_1.$ Here $H_1$ and $H_2$ have the
same elementary theories and, in particular, their
quasi-identities coincide.
\end{thm}

\begin{proof}
Since $H$ is not logically noetherian, we have the inequality
$$LSC(H)\ne q\Var(H).$$
According to \cite{14}, we have a presentation
$q\Var(H)=SCC_{up}(H).$ Here $C_{up}$ is an operator which takes
ultraproducts of algebras.

The class $C_{up}(H_1)$ has an algebra $H_2$ which does not belong to the class $LSC(H_1).$
Therefore, $LSC(H_1)\ne LSC(H_2),$ and the algebras $H_1$ and $H_2$ are not geometrically
equivalent.
The algebra $H_2$ is an ultrapower of the algebra $H=H_1.$
\end{proof}

\begin{thm}\label{thm3}
If the algebras $H_1$ and $H_2$ are logically noetherian, then
they are geometrically equivalent if and only if they have the
same quasi-identities.
\end{thm}

\subsection{Examples. Problems}

Consider the problem  of existing not logically noetherian
algebras in $\Th.$

{\bf 1.} $\Theta=\Grp.$ Let us do it  for the cases $\Th=\Grp$ and
$\Th=\Ass$-$P$, starting with groups. Using \cite{12}, consider
finitely presented groups in the form $F(X)/U,$ where $F(X)$ is a
free group over finite $X,$ and $U$ is an invariant subgroup in
$F=F(X)$ with the finite set of generators (as the invariant
subgroup).

Let $H$ be a discrete direct product of all such $F(X)/U.$ In the
countable group $H$ there is a countable set of finitely generated
subgroups.

Show that the group $H$ is not logically noetherian. We use here
the known Theorem (see \cite{19}) that there exists a continuum of
two-generated simple groups. One of such groups, say $\Gamma,$ is
not embeddable into the group $H$.

Consider a surjection $\mu: F(x,y)\to\Gamma.$
Let $U=\Ker\mu.$
Take a sequence $u_1,u_2,\dots,u_n,\dots$ of all elements of $U.$

Denote by $U_n$ an invariant subgroup in $F=F(X,Y),$ generated by elements $u_1,\dots, u_n.$
The union of all $U_n$ is $U.$
Besides, $F(x,y)/U_n$ is embedded injectively in $H$ and, hence, all $U_n$ are $H$-closed.
We check that $U$ is not an $H$-closed invariant subgroup.

Assume that $U$ is $H$-closed and $\Gamma\approx F(x,y)/U$ is
embedded into $H^I$ for some $I.$ We assume that $\Gamma$ is a
subgroup in $H^I.$ Consider a system of invariant subgroups
$U_\alpha$ in $H^I$ with $H^I/U_\alpha\approx H$ and $\capl_\alpha
U_\alpha=1.$ The image $\Gamma$ in $H^I/U_\alpha$ is isomorphic to
$\Gamma/\Gamma\cap U_\alpha.$ \ $\Gamma $ is a simple group.
 If $\Gamma\cap U_\alpha=\Gamma$ always holds true, then we get a contradiction.
Therefore, $\Gamma\cap U_\alpha=1$ for some $\alpha$ and $\Gamma$ is embedded
into $H\approx
H^I/U_\alpha,$ which contradicts the choice of $\Gamma.$
Thus, the invariant subgroup $U$ is not
$H$-closed, and the group $H$ is not logically noetherian.

Similar considerations are valid in the case $\Th=\Ass$-$P$.

 {\bf 2.} $\Theta=\Ass-P$. Let us call an algebra $H\in\Th$ correct if there exists a surjection
$H\to P.$ A simple algebra is correct if and only if it coincides
with $P.$ In the general case $P$ is a subalgebra in $H.$

Let us note that if $H_\alpha,$\ $\alpha\in I,$ is a family of
correct algebras and
 $H$ is
their free product in $\Th$ then all embeddings $i_\alpha: H_\alpha\to H$ are injective.

Indeed, let us fix $\alpha$ and consider homomorphisms $\nu_\beta:
H_\beta\to H_\alpha$ for all $\beta\in I.$ If $\beta=\alpha,$ then
$\nu_\beta=\nu_\alpha: H_\alpha \to H_\alpha$ is an identical
isomorphism; if $\beta\ne\alpha,$ then $\nu_\beta: H_\beta\to
H_\alpha$ is a homomorphism on the subalgebra $P$ in $H_\alpha.$
By the definition of a free product there is $\nu: H\to H_\alpha$
such that $\nu i_\alpha=\nu_\alpha.$

If now $i_\alpha$ is not injective then we come to the contradiction with the definition of
$\nu_\alpha.$

Consider free algebras  $W=W(X)$ in $\Th=\Ass$-$P$  with finite $X.$
In every such algebra consider finitely-generated ideals $U,$ for which the factor-algebra $W/U$
is correct.
Take for $H$ a free product of all such $W/U$ for different $X$ and $U.$

Assume further that the field $P$ is countable.
Then the algebra $H$ is also countable and in $H$ there exists a countable set of finitely
generated subalgebras.

\begin{thm} The algebra $H$ is not a logically noetherian algebra.
\end{thm}

\begin{proof}.
 We want to show that for the given algebra $H$ there
exists a free algebra $W=W(X)$ such that the union of the
increasing sequence of $H$-closed ideals can be not an $H$ -
closed ideal.

We plan to show this with the help of the appropriate finitely
generated algebra $\Gamma$, which is not embeddable to any
Cartesian power of the algebra $H$. In order to find such an
algebra we need some observations concerning group algebras of
simple groups.

Let $G$ be a simple group and $PG$ be its group algebra. Consider
ideals $V\subset PG$. The canonical homomorphism $PG\to PG/V$
implies $\mu: G\to PG/V$. The kernel $\Ker \mu$ consists of the
elements $g\in G$ such that  $g-1\in V$. This kernel either the
whole group $G$ or $1$. In the first case we have $g-1$ lies in
$V$ for every $g\in G.$  Then the ideal $V$ coincides with the
augmentation ideal $\Delta$. In the second case $(g-1) \in V$
implies $g=1$. In the second case we call the ideal $V$ faithful
ideal.

Thus, if $V\neq \Delta$ then V is faithful. The union of the
increasing sequence of faithful ideals is again a faithful ideal.
Therefore, there are maximal faithful ideals in $PG$. Let $V_O$ be
a maximal faithful ideal. If $V_0 \subset V$ and $V$ does not
coincide with $PG$ then either $V_0=V$ or $V= \Delta$.

Take an algebra $\Gamma=PG/V_0$. We have an injection $G\to
\Gamma$. There are two possibilities:

1. $ \Gamma$ is a simple algebra. 2. $ \Gamma$ has a unique proper
ideal $\bar\Delta=\Delta/V_0$.

Let further the group $G$ be  finitely generated. Then this group
is simultaneously  a finitely generated as a semigroup. Then the
group algebra $PG $ is a finitely generated algebra, and $\Gamma$
is also finitely generated.

The algebra $H$ is countable. In such an algebra there exists not
more than countable set of finitely generated groups. So, we can
find a finitely generated simple group $G$ which is not embeddable
to $H$. Then  the algebra $\Gamma$ is not embeddable to $H$.

Assume that for $\Gamma$ the second case takes place,  that is
there exists a unique ideal $\bar \Delta$ in $\Gamma$. Suppose
that $\Gamma$ is embeddable as a subalgebra to $H^I$. Take a
system of ideals $V_\alpha$, $\alpha\in I$ such that $H^I/V_\alpha
\approx H$ and $\cap_{\alpha \in I}V_\alpha=0$. If $\Gamma\cap
V_\alpha=0$ for some $\alpha$ then $\Gamma$ is embeddable to $H$.
Contradiction. If this intersection is not equal to zero, then it
always contains $\bar \Delta$, which contradicts $\cap_{\alpha \in
I}V_\alpha=0$. Therefore, the algebra $\Gamma$ is not embeddable
in a Cartesian power of $H$.

Take now a finitely generated algebra $W(X)$ with the surjection
$\mu: W(X)\to \Gamma$. Take $U=\Ker \mu$. Let $u_1, \cdots, u_n,
\cdots$ be all elements of algebra $U$. Denote by $U_n$ the ideal
generated by the first $n$ elements. Then algebra $W(X)/U_n$ is
finitely presented and correct. Such algebra is embeddable to $H$.
Hence, every ideal $U_n$ is $H$-closed. However, the union of
these ideals, i.e., the ideal $U$ is not $H$-closed since $W(X)/U$
is isomorphic to $\Gamma$ which is not embeddable to any Cartesian
power of $H$. This means that in the second case we found an
appropriate algebra $W(X)$ which makes $H$ not logically
noetherian.

Suppose now that for $G$ and $\Gamma$ the first case holds, i.e.,
algebra $\Gamma$ is a simple algebra. Consider
$P\times\Gamma=\Gamma^*$.
 There are only two ideals in this algebra, namely $P$ and $\Gamma$.
Besides, assume that $G$ is not embeddable also in $H\times H$.
Since $G$ is embeddable in $\Gamma^*$ then $\Gamma^*$ is not
embeddable in $H$ and $H\times H$. We show that $\Gamma^*$ is not
embeddable in any $H^I$.

As before, consider a system of ideals $V_\alpha$ with
$H^I/V_\alpha\approx H$ and $\cap V_\alpha=0$ in $H^I.$ The image
$\Gamma^*$ in $H^I/V_\alpha$ is isomorphic to the algebra
$\Gamma^*/\Gamma^*\cap V_\alpha.$ There are only two proper ideals
in the algebra $\Gamma^*$, namely, $P$ and $\Gamma.$ If
$\Gamma^*\cap V_\alpha=0$ for some $\alpha, $ then $\Gamma^*$ is
embedded into $H.$ Contradiction. If $\Gamma^*\cap V_\alpha=P$ for
some $\alpha, $ then already $\Gamma $ is embedded into $H.$
Besides, $P\subset H.$ Then $\Gamma^*=P\times\Gamma\subset H\times
H.$ Contradiction. Other cases lead to contradiction with $\cap
V_\alpha=0.$ Hence, the embedding $\Gamma^*\to H^I$ is impossible,

The algebra $\Gamma^*$ is finitely generated and correct. Consider
the surjection $W(X)\to \Gamma^*$. It remains to repeat the
arguments above. We found again the appropriate $W(X)$. So, the
algebra $H$ is not logically noetherian. The theorem is proved.
\end{proof}

The condition on the field to be countable can be eliminated.
Indeed, let $K$ be an arbitrary field and $ P$ its countable
subfield. According to Theorem 4 construct an algebra $H$ over
$P$. It can be proved that by extending scalars to the elements of
$K$ we get an algebra over $K$ which is not logically noetherian.

\begin{prob}\label{prob1}
Let $W=W(X)$ be a free  in $\Ass$-$P$ algebra with $|X|\ge 2.$ Is
it true that $W$ is not geometrically noetherian, but is logically
noetherian?
\end{prob}

\begin{prob}\label{prob2}
The same question for free Lie algebras (see also \cite{26}).
\end{prob}

 \begin{prob}\label{prob3}
Let $G$ be a group, and $PG=H$ its group algebra. The problem is
to find the relationship between the  geometrical and logical
noehterianity for $G$ and $H$.
\end{prob}

This is, indeed, a wide topic related to various problems in the
group algebra theory.

\begin{prob}\label{prob4}
Is it true that there exists continuum  finitely generated simple
associative algebras?\footnote{When the paper was finished I have
been informed that a solution of this problem is contained in the
forthcoming paper [20].}
\end{prob}

\begin{prob}\label{prob5}
Is it true that there exists continuum  finitely generated simple
Lie algebras?
\end{prob}

The positive answer on problem 5 would allow to construct an
example of not logically noetherian Lie algebra.

Note that in the paper by V.Bludov and D.Gusev [8] there is an
example of the solvable group of class 3 which is not logically
noetherian. See also [6], [20],[38].


\section{Geometrical Similarity of Algebras}\label{Geometrical Similarity}

\subsection{Isomorphism of  functors}

If algebras $H_1$ and $H_2$ are geometrically equivalent,
then the categories $C_\Th(H_1)$ and
$C_\Th(H_2)$ coincide, while categories $K_\Th(H_1)$ and $K_\Th(H_2)$ are isomorphic.
The notion of geometrical similarity of algebras is related to necessary and sufficient
conditions of isomorphism of categories of algebraic sets.

Let us recall
the notions of homomorphism and isomorphism of two functors of a
category, which we
will use in the sequel.

Let two functors $\varphi_1,\varphi_2: C_1\to C_2$ of the
categories $C_1,C_2$ be given. The homomorphism (natural
transformation) of functors $s:\varphi_1\to \varphi_2$ is a
function, relating a morphism in $C_2,$ denoted by $s_A:
\varphi_1(A)\to \varphi_2(A)$ to every object $A$ of the category
$C_1.$ For every $\nu: A\to B$ in $C_1$ there is a commutative
diagram

$$
\CD
\varphi_1(A) @> s_A>> \varphi_2(A)\\
@V\varphi_1(\nu)VV @VV\varphi_2(\nu) V\\
\varphi_1(B) @>s_B>> \varphi_2(B)
\endCD
$$
\noindent


in the case of covariant $\varphi_1$ and $\varphi_2.$

For contravariant $\vp_1$ and $\vp_2$ the corresponding diagram is
$$\CD
\varphi_1(B) @> s_B>> \varphi_2(B)\\
@V\varphi_1(\nu) VV @VV\varphi_2(\nu) V\\
\varphi_1(A) @>s_A>> \vp_2(A)
\endCD
$$
\noindent
 An invertible $s:\vp_1\to\vp_2$ is {\it isomorphism
(natural isomorphism) of functors}. The isomorphism property holds
if $s_A:\vp_1(A)\to\vp_2(A)$ is an isomorphism in $C_2$ for any
$A$.

\subsection{Functor $\Cl_H$}

Consider a (contravariant) functor $\Cl_H:\Th^0\to \Set$ for every
algebra $H\in\Th.$ If $W=W(X)$ is an object of $\Th^0,$ then
$\Cl_H(W)$ is the set of all $H$-closed congruences $T$ in $W.$
If, further, $s:W(Y)\to W(X)$ is a morphism of $\Th^0,$ then we
have a mapping of sets $\Cl_H(s):\Cl_H(W(X))\to\Cl_H(W(Y)).$ This
mapping is defined by the following rule: if $T$ is an $H$-closed
congruence in $W(X)$, then $\Cl_H(s)(T)=s^{-1}T.$ It is always an
$H$-closed congruence in $W(Y).$ Here, $w(s^{-1}T)w'$ if
$w^sTw'{}^s.$ For every subvariety $\Th_1$ in $\Th,$ containing an
algebra $H,$ we have also $\Cl_H:\Th_1^0\to\Set.$ These two
different $\Cl_H$ are well correlated. If $W=W(X)\in\Ob\Th^0,$
then $W_0=W_0(X)$ is an object in $\Th_1^0$ with the natural
homomorphism $W\to W_0.$ It is easily checked that there is a
bijection between the sets $\Cl_H(W)$ and $\Cl_H(W_0).$

It follows from definitions that the algebras $H_1$ and $H_2$ are
geometrically equivalent if and only if the functors $\Cl_{H_1}$
and $\Cl_{H_2}$ coincide. Besides, $\Var(H_1)=\Var(H_2).$  The
notion of  {\it geometrical similarity} assumes that the varieties
$\Th_1=\Var(H_1)$ and $\Th_2=\Var(H_2)$ do not necessarily
coincide, but there is an isomorphism of categories $\vp:
\Th_1^0\to\Th_2^0$ with the isomorphism of functors $\alpha(\vp):
\Cl_{H_1}\to\Cl_{H_2}\vp$ depending on $\vp$ under one additional
condition, described later. In the commutative diagram

$$
\CD
\Th_1^0 @>\vp>> \Th_2^0\\
@. @/SE/ \Cl_{H_1}  // @VV \Cl_{H_2} V\\
@. \Set \\
\endCD
$$
\noindent functors $\Cl_{H_1}$ and $\Cl_{H_2}$ act on the
categories $\Th_1^0$ and $\Th_2^0$, respectively,  and
commutativity of the diagram is treated as an isomorphism of
functors $\Cl_{H_1}$ and $\Cl_{H_2}\vp$.

\subsection{Function $\beta$ on the category of free algebras}

Given $\Th$ and the category $\Th^0,$ consider a special function $\beta.$
Take two arbitrary objects $W_1$ and $W_2$ in the category $\Th^0.$
Let $T$ be a congruence in $W_2.$
Denote by $\beta=\beta_{W_1,W_2}(T)$  a binary relation in $\Hom(W_1,W_2).$
This relation is defined as follows:
$s_1\beta s_2$ holds for $s_1,s_2: W_1\to W_2$ if and only if $w^{s_1}Tw^{s_2}$ for every $w\in
W_1.$
The isomorphism $\alpha=\alpha(\vp)$ should commute with the functor $\beta.$

\subsection {Geometrical similarity}

\begin{defn}
Let $H_1$ and $H_2$ be algebras in $\Th,$\ $\Th_1=\Var(H_1),$\
$\Th_2=\Var(H_2).$ The algebras $H_1$ and $H_2$ are called
geometrically similar if
\begin{enumerate}
\item[1.] There exists an isomorphism $\vp: \Th_1^0\to\Th_2^0.$
\item[2.] There exists a function $\alpha=\alpha(\vp)$ such that a bijection
$\alpha(\vp)_W:$

$\Cl_{H_1}(W)\to\Cl_{H_2}(\vp(W))$ holds for every $W\in \Ob
\Th_1^0.$
\item[3.] The function $\alpha$ is coordinated with the function $\beta.$
\end{enumerate}
The last condition means that
$$\vp(\beta_{W_1,W_2}(T))=\beta_{\vp(W_1),\vp(W_2)}(\alpha(\vp)_{W_2}(T)).$$
Here $W_1,W_2$ are objects in $\Th_1^0,$ \ $T$
is an $H_1$-closed congruence in $W_2,$ and for
every relation
$\rho$ in $\Hom(W_1,W_2)$ the relation $\vp(\rho)$ is defined by the rule:
$s_1'\vp(\rho)s_2'$ holds
for $s_1',s_2':\vp(W_1)\to \vp(W_2)$ if there are $s_1,s_2: W_1\to W_2$
such that $\vp(s_1)=s_1',$\ $\vp(s_2)=s_2'$ and $s_1\rho s_2.$

Show now that the function $\alpha=\alpha(\vp)$ is uniquely
determined by these conditions, and give the formula for its
calculation. To this end, consider a function $\rho$ with
$\rho_W=\beta_{W,W}$ for every $W\in\Ob\Th^0.$ Besides, define a
function $\tau$, such that $\tau_W$ is applied to the relation
$\rho$ in $\End W$ for every $W.$ Here $\tau_W(\rho)=T$ is a
relation in $W,$ defined by the rule: $w_1 Tw_2$ if there is $w\in
W$ with $w^\nu=w_1,$\ $w^{\nu'}=w_2$ and $\nu\rho\nu'.$
\end{defn}

It is proved (see \cite{31}) that if $T$ is a congruence in $W,$ then $\tau_W\rho_W(T)=T.$

\begin{prop}\label{prop8}
$\alpha(\vp)_W(T)=\tau_{\vp(W)}\vp(\rho_W(T))$ holds true.
\end{prop}

\begin{proof}
By the condition of coordination between $\alpha$ and $\beta$, we have
$\vp(\rho_W(T))=\rho_{\vp(W)}(\alpha(\vp)_W(T)).$
Here $T$ is an $H_1$-closed  congruence in $W,$ $\alpha(\vp)_W(T)$ is an $H_2$-closed congruence
in $\vp(W).$
Let us apply $\tau_{\vp(W)}.$
$$\tau_{\vp(W)}(\vp(\rho_W(T))=\tau_{\vp(W)}\rho_{\vp(W)}(\alpha(\vp)_W(T))=\alpha(\vp)_W(T).$$
The proved formula allows to state that for every $W\in\Ob \Th_1^0$ the mapping
$\alpha(\vp)_W:\Cl_{H_1}(W)\to\Cl_{H_2}(\vp(W))$ determines the isomorphism of latices of
algebraic sets in $\Hom(W,H_1)$ and $\Hom(\vp(W),H_2).$

Finally, it is proved \cite{32} that $\alpha(\vp)$ gives an
isomorphism of functors $\Cl_{H_1}\to\Cl_{H_2}\vp.$ It is easy to
understand that the relation of geometrical similarity of algebras
in $\Th$ is reflexive, symmetric and transitive. It is also clear
that geometric equivalence is a particular case of geometrical
similarity.
\end{proof}

\subsection{Inner automorphisms of the category of free algebras}

Assume further that $\Var(H_1)=\Var(H_2)=\Th$   for the algebras
$H_1$ and $H_2.$ This is a natural condition; it always holds for
the variety $\Com$-$P$ if the field $P$ is infinite. In this case
similarity of the algebras $H_1$ and $H_2$ is defined by an
automorphism $\vp:\Th^0\to\Th^0.$ For various special $\vp,$
similarity to some extent is reduced to geometrical equivalence.
Let us consider one of the such cases.

For  an arbitrary category $C$, let us call its automorphism $\vp:
C\to C$ an {\it inner}  if there is an isomorphism of functors $s:
1_C\to\vp.$ Here for every object $A$ we have an isomorphism
$s_A:A\to\vp(A)$ and for every $\nu: A\to B$ the diagram
$$\CD
A @> s_A>> \varphi(A)\\
@V\nu VV @VV\varphi(\nu)V\\
B @>s_B>> \vp(B)
\endCD
$$
\noindent
 is commutative. Thus, $\vp(\nu)=s_B\nu s_A^{-1}.$ This
motivates the word ``inner''.

Similarly, one can define an inner
 endomorphism (endofunctor) of a category: it is an arbitrary
$\vp: C\to C,$ isomorphic to a unit automorphism $1_C.$

\begin{prop}\label{prop9}
If similarity of the algebras $H_1$ and $H_2$ is determined by an
inner automorphism $\vp$ of the category $\Th^0,$ then $H_1$ and
$H_2$ are geometrically equivalent.
\end{prop}

\begin{proof}
Let an isomorphism  $s:1_{\Th^0}\to\vp$, an object $W$ in $\Th^0$,
and a congruence $T$ in $W$ be given. Check that
$$\alpha(\vp)_W(T)=s_WT.$$
Here $s_W:W\to\vp(W)$ is an isomorphism of objects and $s_WT$ is a
congruence in $\vp(W)$, defined by the rule: $w_1'(s_WT)w_2'$
holds if and only if $w_1'=s_W(w_1),$\ $w_2'=s_W(w_2)$ and
$w_1Tw_2.$ Denote $s_WT=T^\ast$ and check that
\begin{equation}\label{**}
 \vp(\rho_W(T))=\rho_{\vp(W)}(T^\ast).\tag{$***$}
\end{equation}
Let $\mu,\mu'\in\End\vp(W)$ and $\mu\vp(\rho_W(T))\mu'$ take
place. Then: $\nu,\nu'\in\End W,$\ $\nu\rho_W(T)\nu',$\
$\mu=\vp(\nu),$ \ $\mu'=\vp(\nu').$ For every $w\in W$ we have
$\nu(w)T\nu'(w).$ Thus, $$s_W\nu(w)T^\ast s_W\nu'(w);$$
$$s_W\nu s_W^{-1}(s_Ww)T^\ast s_W\nu's_W^{-1}(s_W(w));$$
$$\mu(w_1)T^\ast\mu'(w_1).$$
Here, $w_1=s_W(w)=w_1$ is an arbitrary element in $\vp(W)$ which
gives us $\mu(\rho_{\vp(W)}(T^\ast))\mu'.$

Let, now, $\mu(\rho_{\vp(W)}(T^\ast))\mu'$ holds. This means that
$\mu(w_1)T^\ast\mu'(w_1)$ holds for every $w_1\in\vp(W).$ Take
$\nu$ and $\nu'\in\End W$ with $\mu=\vp(\nu)=s_W\nu s_W^{-1},$ \
$\mu'=s_W\nu's_W^{-1}$ and $w$ with $s_W(w)=w_1,$ where $w$ is an
arbitrary element in $W.$ We have $s_W(\nu(w))T^\ast
s_W(\nu'(w)).$ This gives $\nu(w)T\nu'(w),$\ $\nu\rho_W(T)\nu'.$
Then $\vp(\nu)\vp(\rho_W(T))\vp(\nu')$ and
$\mu\vp(\rho_W(T))\mu'.$ The equality \eqref{**} is checked. Now
we have
$$\alpha(\vp_W(T))=\tau_{\vp(W)}\vp(\rho_W(T))
=\tau_{\vp(W)}\rho_{\vp(W)}(T^\ast)=T^\ast=s_WT.$$

Let $T$ be an $H_1$-closed congruence in $W.$
Then $s_WT$ is an $H_2$-closed congruence in
$\vp(W)$ by the definition of similarity.
On the other hand,
using the isomorphism $s_W^{-1}:\vp(W)\to W$ and the fact that isomorphism of
objects in $\Th^0$ preserves the
$H$-closeness condition for every $H$ (see \cite{29}), we
conclude that $T$ is an $H_2$-closed congruence as well.
Hence every $H_1$-closed congruence in $W$ is $H_2$-closed.
Applying $\vp^{-1},$ conclude the  opposite.
Thus, $H_1$ and $H_2$ are geometrically equivalent.
\end{proof}

Other examples of this kind will be
given in the section devoted to the case $\Th=\Ass$-$P$.
It is proved for $\Th=\Grp$
that all automorphisms of the category $\Th^0$ are inner \cite{23}.

\section{Geometrical Coordination of Algebras}\label{geocoor}

\subsection{Additional information on categories}

Coordination of algebras leads to the necessary and sufficient conditions of equivalence of two
categories of algebraic sets.

We recall here some required information from category theory [25], [39].

Let the categories $C_1$ and $C_2$ be given. They are {\it
equivalent} if there exists a pair of functors $\vp: C_1\to C_2$
and $\psi: C_2\to C_1$ such that $\psi\vp\approx 1_{C_1},$\
$\vp\psi\approx 1_{C_2}.$ Here, $1_C$ is a unity functor of a
category. The sign $\approx$ here denotes isomorphism of functors.
We say that the pair $(\vp,\psi)$ determines {\it equivalence of
categories}  $C_1$ and $C_2.$ If $\psi\vp=1_{C_1},$\
$\vp\psi=1_{C_2},$  then the pair $(\vp,\psi)$ determines {\it
isomorphism of categories} and $\psi=\vp^{-1}.$

It is proved \cite{39} that if $(\vp,\psi)$ is an equivalence, then
 each of the  functors $\vp$ and $\psi$ possesses  the
following two properties

1. \ Completeness;

2.\ Univalencity.

For $\vp: C_1\to C_2$, completeness means that for every object
$B$ of $C_2$ there exists an object $A$ of $C_1,$ such that
$\vp(A)\approx B.$ Univalencity means that for any two objects $A$
and $B$ of $C_1$, the functor $\vp$ induces a bijection
$\vp_{A,B}: Hom(A,B)\to\Hom(\vp(A),\vp(B)).$ In particular, for
every $A$ this gives an isomorphism $\vp_A: \End A\to \End\vp(A).$
Let us call a functor $\vp$ with these two properties {\it a
relational isomorphism} of categories. If $\vp$ is a relational
isomorphism, then it has a relational inverse functor $\psi$ such
that the pair $(\vp,\psi)$ determines equivalence of categories.
There could be many relational inverse functors for $\vp.$

If $\vp$ is an isomorphism, then there is only one inverse functor $\vp^{-1},$ but there are many
relational inverse ones.

Let us fix  a small category $C.$ Consider endofunctors
(endomorphisms) $\vp: C\to C.$ They constitute a semigroup $\End
C.$ Relational automorphisms (autoequivalences) form a
subsemigroup in $\End C,$ denoted by $\widetilde{\Aut}(C).$ The
group of automorphisms $\Aut(C)$ is a group of invertible elements
in $\End(C).$

It is checked that the isomorphism relation $\approx$ in is a
congruence in  $\End(C)$: $\vp_1\approx \vp_2$ and
$\psi_1\approx\psi_2$ imply $\vp_1\psi_1\approx \vp_2\psi_2.$

We can now pass to the factor-semigroup
$\End^0(C)=\End(C)/\approx.$ Denote a group of invertible elements
in $\End^0(C)$ by $\Aut^0(C).$ Fix a natural homomorphism
$\dl:\End(C)\to\End^0(C).$ If $\dl(\vp)=\ov\vp$ is an invertible
element in $\Aut^0(C),$ then we take $\ov\psi=\ov\vp^{-1}$ and,
therefore, $\ov\vp\ov\psi=\ov\psi\ov\vp=\ov{1_C},$\ $\ov{\vp\psi}
=\ov{\psi\vp}=\ov{1_C},$\ $\vp\psi\approx 1_C\approx\psi\vp.$
Thus, $\vp$ is an autoequivalence, like $\psi.$ The subsemigroup
$\widetilde{\Aut}(C)$ is a full co-image of the group $\Aut^0(C).$
A homomorphism $\dl$ induces the homomorphism $\dl:\Aut(C)\to
\Aut^0(C).$ The latter is surjective if every autoequivalence
$\vp$ of the category $C$ is isomorphic  to some automorphism
$\psi,$\ $\vp\approx\psi,$\ $\vp\psi^{-1} =\vp_0\approx 1_C,$\
$\vp=\vp_0\psi.$ Here $\vp_0$ is also an autoequivalence , namely
it is an inner one. Besides, let us note that if $\vp=\vp_0\psi,$
then all relationally inverse functors to $\vp$ are of the form
$\psi^{-1}\vp_1,$ where $\vp_1$ is an arbitrary functor,
isomorphic to $1_C.$

Every $\vp_1, $ isomorphic to $1_C$ is simultaneously an
autoequivalence of the category.

Note that the kernel of
the homomorphism $\dl:\Aut(C)\to\Aut^0(C)$ is an invariant subgroup
$\Int(C)$ in $\Aut(C),$ consisting of all inner automorphisms, which is
 isomorphic to a trivial automorphism.

Given a small category $C$ and an object $A,$  denote by $[A]$ the
class (set) of all objects in $C,$ isomorphic to $A.$ The set of
all objects $\Ob(C)$ is decomposed into such classes.

\begin{thm}\label{thm4} (G. Zhitomirsky  \cite{40})
If all classes $[A]$ have the same cardinality pairwise, then every autoequivalence of the
category $C$ is isomorphic to an automorphism.
\end{thm}

\begin{proof}
Let $\vp: C\to C$ be an autoequivalence.
For every object $A$ we set:  $\ov\vp[A]=[\vp(A)].$

It follows from the general categorical considerations that
$\ov\vp$ is a substitution on the set of classes of isomorphic
objects: its definition does not depend on the choice of the
representative $A$ in the classes of isomorphic objects. In the
conditions of the theorem we have a bijection
$\psi_{[A]}:[A]\to\ov\vp[A].$ Fix these bijections. Further, for
every object $A$ we set:
$$\psi(A)=\psi_{[A]}(A)\in\ov\vp[A]=[\vp(A)].$$

Here $\psi$ is a substitution on the set $\Ob C.$
Since $\psi(A)\in[\vp(A)],$ then $\psi(A)$ and $\vp(A)$ are isomorphic.
For every $A$ fix some isomorphism $s_A:\vp(A)\to\psi(A).$
For any $\nu: A\to B$ consider a diagram
$$
\CD
\varphi (A) @> s_A>> \psi(A)\\
@V\varphi(\nu)VV @VV\psi(\nu)V\\
\varphi (B) @>s_B>> \psi(B)
\endCD
$$
\noindent
 Correspondingly, $\psi(\nu)=s_B\vp(\nu)s^{-1}_A.$ Under
such a definition, $\psi$ is an automorphism of the category $C,$
the diagram is commutative and $\vp$ and $\psi$ are isomorphic.
\end{proof}

Let us apply these general facts to the category $\Th^0,$ where
$\Th$ is an arbitrary variety of algebras. Here for every algebra
$W=W(X)$, the class $[W]$ has the same cardinality as the initial
universal set $X^0.$ Hence, we have the following corollary.

\begin{cor}
Every autoequivalence of the category $\Th^0$ is isomorphic to an
 automorphism.
\end{cor}

For every autoequivalence  $\vp$ we have an automorphism $\psi$
with $\vp=\vp_0\psi,$ where $\vp_0$ is an inner autoequivalence.
The homomorphism $\dl:\Aut(\Th^0)\to\Aut^0(\Th^0)$ is always
surjective.

It is easy to understand that \thmref{thm4} admits a
generalization. Every equivalence of different categories
$\Th_1^0$ and $\Th_2^0$ is naturally isomorphic to an isomorphism
of these categories. If $\Th_1^0$ and $\Th_2^0$ are equivalent,
then they are isomorphic.

Fix an arbitrary object $A_0$ in every class $[A].$ This gives the
full subcategory in $C.$ Such a subcategory is considered as the
skeleton of the category $C,$ denoted by $\tilde C.$ The category
$\tilde C$ can be represented also as a category of classes $[A].$

An autoequivalence $\vp: C\to C$ is called special if
$\ov\vp[A]=[\vp(A)]=[A]$ for any $A.$
This means that the objects $A$ and $\vp(A)$ are always isomorphic.

\begin{thm}\label{thm5}
If $\vp$ is a special autoequivalence, then it can be represented
as $\vp=\vp_0\vp_1,$ where $\vp_0$ is an inner autoequivalence and
$\vp_1$ is an automorphism which does not change objects.
\end{thm}

\begin{proof}
We do not use here the previous theorem and the axiom of choice.
First we build an inner autoequivalence $\vp_0,$ setting
$\vp_0(A)=\vp(A)$ for every object $A.$
Then we fix an isomorphism $s_A: A\to\vp(A)=\vp_0(A).$
For $\nu: A\to B$ we set
$$\vp_0(\nu)=s_B\nu s_A^{-1}:\vp_0(A)\to\vp_0(B).$$
Here $\vp_0$ is a functor and $\vp_0\approx 1_C.$
Solving the equation $\vp=\vp_0\vp_1$ with respect to $\vp_1,$
we set $\vp_1(A)=A$ for every $A.$
For $\nu: A\to B$ the equality $\vp(\nu)=\vp_0\vp_1(\nu)$ should hold; here we have $\vp_1(\nu):
A\to B$ and $\vp_0\vp_1(\nu)=s_B\vp_1(\nu)s_A^{-1}.$
Setting $\vp_1(\nu)=s_B^{-1}\vp(\nu)s_A,$ we find the automorphism
$\vp_1,$ which solves the
equation.
\end{proof}

\subsection{Geometrical coordination}

Let us pass to the notion of geometrical coordination of algebras, generalizing geometrical
similarity.

\begin{defn}\label{def6}
Let the algebras $H_1$ and $H_2$
 be given in $\Th,$ \ $\Th_1=\Var(H_1),$\ $\Th_2=\Var(H_2).$
The algebras $H_1$ and $H_2$ are called coordinated if
\begin{enumerate}
\item[1)] There exists an equivalence of categories
$\vp: \Th_1^0\to\Th_2^0$ and
$\psi:\Th_2^0\to\Th_1^0$
\item[2)] For the pair $(\vp,\psi)$ there exist embeddings:
$$\alpha(\vp)_W:\Cl_{H_1}(W)\to\Cl_{H_2}(\vp(W)),\quad
\quad W\in\Ob\Th_1^0;$$
$$\alpha(\psi)_W:\Cl_{H_2}(W)\to\Cl_{H_1}(\psi(W)),\quad\quad
W\in\Ob\Th_2^0.$$
\item[3)] The functions $\alpha(\vp)$ and $\alpha(\psi)$ commute with
the corresponding $\beta.$
\end{enumerate}

It follows from the third condition that, in particular, if $W\in\Ob\Th_1^0$ and $T$ is an
$H_1$-closed congruence in $W,$ then $\varphi(\rho_W(T))=\rho_{\vp(W)}(\alpha(\vp)_W(T)).$
As above, we deduce formulas for the corresponding $W$ and $T:$
$$\alpha(\vp)_W(T)=\tau_{\vp(W)}\vp(\rho_W(T)),$$
$$\alpha(\psi)_W(T)=\tau_{\psi(W)}\psi(\rho_W(T)),$$
\end{defn}

In the proofs that follow, we sometimes take into account the
univalencity property of the functors $\vp$ and $\psi.$

Note that from the definition follows that the transitions
$$\alpha(\vp):\Cl_{H_1}\to\Cl_{H_2}\vp,$$
$$\alpha(\psi):\Cl_{H_2}\to\Cl_{H_1}\psi.$$
turn out to be natural transformations of functors.

\begin{prop}\label{prop10}
If $\vp:\Th_1^0\to\Th_2^0$  is an isomorphism of categories and
$\psi=\vp^{-1}$, then the coordination of the algebras $H_1$ and
$H_2$ means  that these algebras are similar.
\end{prop}

\begin{proof}
We need to check that in the conditions above
$$\alpha(\vp)_W: \Cl_{H_1}(W) \to \Cl_{H_2}(\vp(W))$$
is a bijection, and
$$\alpha(\psi)_{\vp(W)}: \Cl_{H_2}(\vp(W))\to\Cl_{H_1}(W)$$
is the inverse bijection, $W\in\Th_1^0.$

Take $W\in\Ob\Th_1^0$ and an $H_1$-closed congruence $T$ in $W.$
Then
$$\vp(\rho_W(T))=\rho_{\vp(W)}(\alpha(\vp)_W(T)),$$
where $\alpha(\vp)_W(T)$ is $H_2$-closed congruence in $\vp(W)$.
Applying $\psi=\vp^{-1},$ we get
$$\rho_W(T)=\psi(\rho_{\vp(W)}(\alpha(\vp)_W(T))=\rho_W(\alpha(\psi)_{\vp(W)}\alpha(\vp)_W(T)).$$
Hence, $T=\alpha(\psi)_{\vp(W)}\alpha(\vp)_W(T).$

We get a similar result if we take $W\in\Th_2^0$ and $T$ is an $H_2$-closed congruence in $W.$

Evidently, the coordination relation of two algebras is reflexive
and symmetric. Transitivity follows from the considerations below.

Let the algebras $H_1$, $H_2$ and $H_3$  be given in the variety
$\Th.$ Correspondingly, $\Th_1=\Var(H_1),$\ $\Th_2=\Var(H_2),$\
$\Th_3=\Var(H_3).$ Let the pair of functors
$\vp_1:\Th_1^0\to\Th_2^0$ and $\psi_1:\Th_2^0\to\Th_1^0$ determine
coordination of the algebras $H_1$ and $H_2$, and another pair
$\vp_2:\Th_2^0\to\Th_3^0,$\ $\psi_2:\Th_3^0\to\Th_2^0$ determines
coordination for $H_2$ and $H_3.$

We have $\vp=\vp_2\vp_1:\Th_1^0\to\Th_3^0$ and $\psi=\psi_1\psi_2:
\Th_3^0\to\Th_1^0$. Check that the pair $(\vp,\psi)$ determines
coordination of the algebras $H_1$ and $H_2.$

Calculate $\alpha(\vp_2\vp_1)$ and $\alpha(\psi_1\psi_2).$
Take $W\in\Ob\Th_1^0,$ and let $T$ be an $H_1$-closed congruence in $W.$
Let us consider the congruence
$$\alpha(\vp_2\vp_1)_W(T)=\tau_{\vp_2\vp_1(W)}\vp_2\vp_1(\rho_W(T)).$$
We have (compare \propref{prop10}) $\vp_1(\rho_W(T))=\rho_{\vp_1(W)}(\alpha(\vp_1)_W(T)),$ where
$\alpha(\vp_1)_W(T)$ is an $H_2$-closed congruence in $\vp_1(W).$
Further,
$$\vp_2\vp_1(\rho_W(T))=\vp_2(\rho_{\vp_1(W)}(\alpha(\vp_1)_W(T))=\rho_{\vp_2\vp_1(W)}(\alpha(\vp_2)_
{\vp_1(W)}\alpha(\vp_1)_W(T)).$$ Applying $\tau_{\vp_2\vp_1(W)}$
we get $\alpha(\vp_2\vp_1)_W(T)=\alpha(\vp_2)_{\vp_1(W)}\alpha(\vp
_1)_WT.$ Here the congruence $\alpha(\vp_2\vp_1)_W(T)$ is an
$H_3$-closed congruence in $\vp_2\vp_1(W),$ since
$\alpha(\vp_1)_W(T)$ is an $H_2$-closed congruence in $\vp_1(W).$
We have an inclusion
$$\alpha(\vp_2\vp_1)_W:\Cl_{H_1}(W)\to\Cl_{H_3}(\vp_2\vp_1(W)).$$
Similarly, we calculate $\alpha(\psi_1\psi_2)_W(T)$ for
$W\in\Ob\Th_3^0,$ where $T$ is an $H_3$-closed congruence in $W.$
This gives an embedding
$$\alpha(\psi_1\psi_2)_W:\Cl_{H_3}(W)\to\Cl_{H_3}(\psi_1\psi_2(W)).$$
Commutativity of $\alpha$ and $\beta$ is evident.
This gives the corresponding transitivity.
\end{proof}

\begin{prop}\label{prop11}
Let $\Var(H_1)=\Var(H_2)=\Th,$ and $(\vp, \psi)$ be an
autoequivalence of the category $\Th^0,$ and $H_1$ and $H_2$ be
coordinated algebras in respect to $(\phi, \psi).$ If this
autoequivalence is inner, then $H_1$ and $H_2$ are geometrically
equivalent.
\end{prop}

\begin{proof}
The proof is similar to that for geometrical similarity. We take
into account univalencity of the functors $\vp$ and $\psi.$
\end{proof}

\subsection{Decomposition of similarity and coordination relations}

Let the pair of functors $(\vp,\psi)$ determine coordination of
the algebras $H_1$ and $H_2,$ and there is a decomposition
$\vp=\vp_0\vp_1,$\ $\psi=\psi_1\psi_0.$ Assume also that there is
an algebra $H$ such that the pair $(\vp_1,\psi_1)$ determines
coordination of the algebras $H$ and $H_1.$ We want the pair
$(\vp_0,\psi_0)$ to determine coordination of the algebras $H$ and
$H_2.$

We solve this problem of decomposition of coordination relations
in the conditions:
\begin{enumerate}
\item[1)] $\Var(H_1)=\Var(H)=\Var(H_2)=\Th.$

\item[2)] $\vp_1=\zeta,$\ $\psi_1=\zeta^{-1},$ where $\zeta$ is an automorphism of the category
$\Th^0,$ determining similarity of the algebras $H_1$ and $H.$

\item[3)] The automorphism $\zeta$ does not change objects.

\end{enumerate}

The next proposition is valid under the conditions (1) -- (3).

\begin{prop}\label{prop12}
Let the pair $(\vp,\psi)$ determine coordination of the algebras
$H_1$ and $H_2$, \ $\vp=\vp_0\zeta,$\ $\psi=\zeta^{-1}\psi_0,$ and
the automorphism $\zeta$ determine similarity of the algebras
$H_1$ and $H.$ Then the pair $(\vp_0,\psi_0)$ is an
autoequivalence of the category $\Th^0,$ determining coordination
of the algebras $H$ and $H_2.$
\end{prop}

\begin{proof}
Taking into account conditions, consider
$$\alpha(\vp)_W=\alpha(\vp_0\zeta)_W=\alpha(\vp_0)_{\zeta(W)}\alpha(\zeta)_W=\alpha(\vp_0)_W\cdot
\alpha(\zeta)_W.$$ We  have an embedding
$\alpha(\vp)_W:\Cl_{H_1}(W)\to \Cl_{H_2}(\vp(W)).$ Besides,
$\vp(W)=\vp_0(W).$ We have also a bijection
$\alpha(\zeta)_W:\Cl_{H_1}(W)\to\Cl_H(W).$

Let now $T^\ast$ be an arbitrary $H$-closed congruence in $W,$
$T^\ast\in\Cl_H(W).$ We take
 $T^\ast=\alpha(\zeta)_W(T),$ \
$T\in\Cl_{H_1}(W).$ Then
$$\alpha(\vp_0)_W(T^\ast)=\alpha(\vp_0)_W\alpha(\zeta)_W(T)=\alpha(\vp)_W(T)\in\Cl_{H_2}(\vp_0(W)).$$
Thus, we have an embedding $\alpha(\vp_0)_W:
\Cl_H(W)\to\Cl_{H_2}(\vp_0(W)).$

We then work with $\psi=\zeta^{-1}\psi_0.$

For every algebra $W\in\Ob\Th^0$ we have an embedding
$\alpha(\psi)_W: \Cl_{H_2}(W)\to\Cl_{H_1}(\psi(W)).$ Here,
$\psi(W)=\psi_0(W),$\ $\psi_0=\zeta\psi.$
Further,
$$\alpha(\psi_0)_W=\alpha(\zeta)_{\psi(W)}\alpha(\psi)_W=\alpha(\zeta)_{\psi_0(W)}\alpha(\psi)_W.$$

Let now $T\in\Cl_{H_2}(W).$
Then $\alpha(\psi)_W(T)\in\Cl_{H_1}(\psi(W))=\Cl_{H_1}(\psi_0(W)).$
We have also a bijection
$\alpha(\zeta)_{\psi_0(W)}:\Cl_{H_1}(\psi_0(W)\to \Cl_H(\psi_0(W)).$
Hence, $\alpha(\psi_0)_W(T)=\alpha(\zeta)_{\psi_0(W)}\alpha(\psi)_W(T)$ and
$\alpha(\psi_0)_W(T)\in\Cl_H(\psi_0(W)).$
This means that there exists an embedding $\alpha(\psi_0)_W:\Cl_{H_2}(W)\to\Cl_H(\psi_0(W)).$

Note further that the pair $(\vp_0,\psi_0)$ is an autoequivalence
of the category $\Th^0.$ We have:
$$\vp\psi=\vp_0\zeta\zeta^{-1}\psi_0=\vp_0\psi_0\approx 1_{\Th^0}$$
$$\psi\vp=\zeta^{-1}\psi_0\vp_0\zeta\approx 1_{\Th^0}\quad\text{and}\quad \psi_0\vp_0\approx
1_{\Th^0}.$$
Besides, we have checked that there are embeddings
$$\alpha(\vp_0)_W:\Cl_H(W)\to \Cl_{H_2}(\vp_0(W)).$$
$$\alpha(\psi_0)_W:\Cl_{H_2}(W)\to\Cl_H(\psi_0(W)).$$

It is left to check coordination of the mappings $\alpha(\vp_0)_W$
and $\alpha(\psi_0)_W$ with the function $\beta.$ Let $W_1$ and
$W_2$ be objects in $\Th^0$ and $T$ a congruence in $W_2.$ Denote
$\beta_{W_1,W_2}(T)=\beta,$ and let $\alpha(\vp)_{W_2}(T)=T^\ast$
be a congruence in $\vp(W_2).$ Denote $
\beta^\ast=\beta_{\vp(W_1),\vp(W_2)}(T^\ast).$

We consider $T$ to be an $H_1$-closed congruence in $W_2.$ Then
$T^\ast$ is an $H_2$-closed congruence in $\vp(W_2).$ Under these
conditions $\vp(\beta)=\beta^\ast.$ We will repeat the similar
calculations  for $\vp_0.$ Proceed from $\vp=\vp_0\zeta$ and
$\alpha(\vp)_{W_2}=\alpha(\vp_0)_{W_2}\alpha(\zeta)_{W_2}.$ Let
now $T$ be an $H$-closed congruence in $W_2,$\
$T=\alpha(\zeta)_{W_2}(T_1)$ where $T_1$ is an $H_1$-closed
congruence in $W_2.$ We have
$\alpha(\vp)_{W_2}(T_1)=\alpha(\vp_0)_{W_2}(T).$ Besides,
$\vp(W_1)=\vp_0(W_1),$\ $\vp(W_2)=\vp_0(W_2).$ Then
$$\beta_{\vp_0(W_1),\vp_0(W_2)}(\alpha(\vp_0)_{W_2}(T))=\beta_{\vp(W_1),\vp(W_2)}(\alpha(\vp)_{W_2}
(T_1))=\vp(\beta_{W_1,W_2}(T_1)).$$ We need to check that
$\vp(\beta_{W_1,W_2}(T_1))=\vp_0(\beta_{W_1,W_2}(T)).$ Here
$T=\alpha(\zeta)_{W_2}(T_1).$ The functor $\zeta$ is coordinated
with $\beta.$ Hence,
$\zeta(\beta_{W_1,W_2}(T_1))=\beta_{W_1,W_2}(\alpha(\zeta)_{W_2}(T_1)),$
and $T_1$ is an $H_1$-closed congruence in $W_2.$ Applying
$\vp_0,$\
$$\vp_0\zeta(\beta_{W_1,W_2}(T_1))=\vp_0(\beta_{W_1,W_2}(\alpha(\zeta)_{W_2}(T_1))=\vp_0(\beta_{W_1,W_2}
(T)).$$ This gives
$\vp(\beta_{W_1,W_2}(T_1))=\vp_0(\beta_{W_1,W_2}(T)).$ Finally,
$\beta_{\vp_0(W_1),\vp_0(W_2)}(\alpha(\vp_0)_{W_2}(T))=\vp_0(\beta_{W_1,W_2}(T)).$
We have checked coordination of $\vp_0$ and $\beta.$

Let us pass to $\psi_0$ and $\beta,$\ $\psi_0=\zeta\psi.$ Use that
the functors $\zeta$ and $\psi$ commute with $\beta.$
 Take once more the objects $W_1$ and $W_2$ in $\Th^0,$ and let $T$ be an $H_2$-closed congruence
in $W_2.$  We have
$\psi(\beta_{W_1,W_2}(T))=\beta_{\psi(W_1),\psi(W_2)}(\alpha(\psi)_{W_2}(T)).$
Applying $\zeta,$ we get
$$\psi_0(\beta_{W_1,W_2}(T))=\zeta(\beta_{\psi(W_1),\psi(W_2)}(\alpha(\psi)_{W_2}(T))=$$
$$=\beta_{\psi_0(W_1),
\psi_0(W_2)}(\alpha(\zeta)_{\psi(W_2)}\alpha(\psi)_{W_2}(T))
=\beta_{\psi_0(W_1),\psi_0(W_2)}(\alpha(\psi_0) _{W_2}(T)).$$ We
have also checked correspondence of $\psi_0$ and $\beta.$ Thus,
$\vp_0$ and $\psi_0$ determine geometrical coordination of the
algebras $H$ and $H_2.$ The proposition is proved.
\end{proof}

\section{Isomorphisms and equivalences of categories of algebraic sets}\label{IsoEq}

\subsection{Correctness}

Define first correct isomorphisms and correct equivalences under the conditions
$\Var(H_1)=\Var(H_2)=\Th.$
Every isomorphism
$$F: K_\Th(H_1)\to K_\Th(H_2)$$
is in one-to-one correspondence with the isomorphism
$$\Phi: C_\Th(H_1)\to C_\Th(H_2).$$
The category of affine spaces $K_\Th^0(H)$ is a subcategory of $K_\Th(H).$

{\it Correctness of an isomorphism }$F$ assumes that $F$ respects
the categories of affine spaces, that is $F$ induces
$$F^0: K_\Th^0(H_1)\to K_\Th^0(H_2).$$
Correspondingly, $\Phi$ induces an automorphism of the category
$\Th^0$
$$\vp: \Th^0\to\Th^0.$$
Recall that the category $\Th^0$ is a subcategory in $C_\Th(H_1)$
and $C_\Th(H_2).$ The equality $F(\Hom(W,H_1))=\Hom(\vp(W),H_2)$
always holds true. Besides, suppose that for every object $(X,A)$
of the category $K_\Th(H_1)$,\ the equality $F((X,A))=(Y,B),$
where $B$ is an algebraic set in the affine space
$\Hom(W(Y),H_2),$ and $W(Y)=\vp(W(X))$ holds.

This definition of correctness of   isomorphism is quite natural
and in the sequel, isomorphism of categories  means correct
isomorphism.

Note that it follows from the definition that if $\mu:W(X)\to
W(X)/T$ is a natural homomorphism in the category $C_\Th(H_1),$
then a natural homomorphism $\Phi(\mu):\vp(W)\to\vp(W)/T^\ast$ in
$C_\Th(H_2)$ corresponds to this $\mu$ [33].

Let us now pass to the correct equivalence.
We have
$$F_1: K_\Th(H_1)\to K_\Th(H_2),$$
$$F_2: K_\Th(H_2)\to K_\Th(H_1).$$
The pair of functors $(F_1,F_2)$ determines the equivalence of categories.
Simultaneously, we have an equivalence
$$\Phi_1: C_\Th(H_1)\to C_\Th(H_2),$$
$$\Phi_2: C_\Th(H_2)\to C_\Th(H_1).$$
As we have done above, we assume correspondence of the functors
with the categories of affine spaces and,  subsequently, with
$\Th^0.$ In particular, the pair $(\Phi_1,\Phi_2)$ induces
autoequivalence of the category $\Th^0.$ The functors $\vp:
\Th^0\to\Th^0$ and $\psi: \Th^0\to\Th^0$ are relatively mutually
inverse. Here, as before, the functors $\Phi_1$ and $\Phi_2$ are
coordinated with natural homomorphisms.

\subsection{Isomorphism and equivalence of categories}

The following two theorems are of universal character; they relate
to arbitrary varieties $\Th.$ Their usage assumes knowledge of the
structure of automorphisms and autoequivalences of categories
$\Th^0$ in various special situations.

\begin{thm}\label{thm6}
The categories $K_\Th(H_1)$ and $K_\Th(H_2)$ are isomorphic if and
only if the algebras $H_1$ and $H_2$ are geometrically similar.
\end{thm}

\begin{proof} See \cite{31}.\end{proof}

\begin{thm}\label{thm7}
The categories $K_\Th(H_1)$ and $K_\Th(H_2)$ are equivalent if and only if the algebras $H_1$ and
$H_2$ are geometrically coordinated.
\end{thm}
\begin{proof}
Let, first, $K_\Th(H_1)$ and $K_\Th(H_2)$ be (correctly)
equivalent. We use the equivalence
$$\Phi=\Phi_1: C_\Th(H_1)\to C_\Th(H_2),$$
$$\Psi=\Phi_2: C_\Th(H_2)\to C_\Th(H_1).$$
This pair induces the autoequivalence
$$\vp:\Th^0\to\Th^0,\quad \psi: \Th^0\to\Th^0,\quad \vp\psi\approx 1_{\Th^0}\approx \psi\vp.$$

Check that the pair $(\vp,\psi)$ determines coordination of the
algebras $H_1$ and $H_2.$ Take functions $\alpha(\vp)$ and
$\alpha(\psi).$ Let $T$ be an $H_1$-closed congruence in $W,$\
$W\in\Ob\Th^0.$ Consider a natural homomorphism $\mu: W\to W/T.$
It is a morphism in $C_\Th(H_1)$ with a corresponding natural
homomorphism $\Phi(\mu): \vp(W)\to \vp(W)/T^\ast.$ The congruence
$T^\ast$ is $H_2$-closed and uniquely defined.

Setting $\alpha(\vp)_W(T)=T^\ast,$ we have
$\alpha(\vp)_W:\Cl_{H_1}(W)\to\Cl_{H_2}(\vp(W)).$ Similarly,
$\alpha(\psi)_W:\Cl_{H_2}(W)\to\Cl_{H_1}(\psi(W)).$

Now  we need to check commutativity with the function $\beta.$
Take $W_1,W_2\in\Ob \Th^0.$ Let $T$ be an $H_1$-closed congruence
in $W_2.$ Consider a natural homomorphism $\mu_T:W_2\to W_2/T.$
For $ s_1,s_2: W_1\to W_2$ the relation $s_1\beta_{W_1,W_2}(T)s_2$
holds if and only if the equality $\mu_Ts_1=\mu_Ts_2$ takes place.
Rewrite in these terms the corresponding commutativity condition.
Given $\mu_Ts_1=\mu_Ts_2,$ apply $\Phi$ to the  equality above and
get $\Phi(\mu_T)\vp(s_1)=\Phi(\mu_T)\vp(s_2).$ Denote
$\Phi(\mu_T)=\mu_{T^\ast}:\vp(W_2)\to \vp(W_2)/T^\ast.$ This is a
natural homomorphism with $T^\ast=\alpha(\vp)_{W_2}(T).$ We have
$\mu_{T^\ast}\vp(s_1)=\mu_{T^\ast}\vp(s_2).$ This is equivalent to
\begin{equation}\label{****}
\vp(s_1)\beta_{\vp(W_1),\vp(W_2)}(T^\ast)\vp(s_2),\tag{$****$}
\end{equation}
  and $s_1\beta_{W_1,W_2}(T)s_2$ implies
  $$
  \varphi(s_1)\beta_{\varphi(W_1),\varphi(W_2)}(\alpha(\varphi)_{W_2}(T))\varphi(s_2).
  $$

Let now $
\mu_{T^\ast}s_1'=\mu_{T^\ast}s_2'$ hold for $s_1',s_2':\vp(W_1)\to\vp(W_2).$
Using univalencity of the functor $\vp,$ find $s_1,s_2: W_1\to W_2$ with $\vp(s_1)=s_1',$\
$\vp(s_2)=s_2'.$
Then
$$\mu_{T^\ast}\vp(s_1)=\mu_{T^\ast}\vp(s_2);$$
$$\Phi(\mu_T)\vp(s_1)=\Phi(\mu_T)\vp(s_2);$$
$$\Phi(\mu_Ts_1)=\Phi_(\mu_Ts_2).$$
Using univalencity of the functor $\Phi,$ we conclude:
$\mu_Ts_1=\mu_Ts_2.$ Hence, the condition
$s_1'\beta_{\vp(W_1),\vp(W_2)}(T^\ast)s_2'$ holds if and only if
$s_1'=\vp(s_1),$\ $s_2'=\vp(s_2)$ and $s_1\beta_{W_1,W_2}(T)s_2.$
This means exactly that
$$\vp(\beta_{W_1,W_2}(T))=\beta_{\vp(W_1),\vp(W_2)}(T^\ast)=\beta_{\vp(W_1),\vp(W_2)}(\alpha(\vp)_{W_2}(T)).$$
The commutativity condition for $\alpha$ and $\beta$ is checked.

The proof for the functor $\psi: \Th^0\to\Th^0$ is similar.

We have proved that the algebras $H_1$ and $H_2$ are geometrically
coordinated.

Prove the opposite.

Let the algebras $H_1$ and $H_2$ be geometrically coordinated.
Prove that the categories $K_\Th(H_1)$ and $K_\Th(H_2)$ are
correctly equivalent. It is sufficient to prove this for the
categories $C_\Th(H_1)$ and $C_\Th(H_2).$

Proceed from the autoequivalence $\vp: \Th^0\to\Th^0,$\ $\psi:\Th^0\to\Th^0$ and the
corresponding functions $\alpha(\vp)$ and $\alpha(\psi).$

For every $W\in\Th^0$ we have the mappings
$$\alpha(\vp)_W: \Cl_{H_1}(W)\to\Cl_{H_2}(\vp(W)),$$
$$\alpha(\psi)_W: \Cl_{H_2}(W)\to\Cl_{H_1}(\psi(W)).$$
Define the functors $\Phi: C_\Th(H_1)\to C_\Th(H_2)$ and $\Psi:
C_\Th(H_2)\to C_\Th(H_1)$ with the conditions $\Psi\Phi\approx
1_{C_\Th(H_1)}$ and $\Phi\Psi\approx 1_{C_\Th(H_2)},$ determining
correct equivalence of the categories $C_{\Th(H_1)}$ and
$C_\Th(H_2).$

Start with the definition of $\Phi.$ Take an arbitrary object
$W/T,$ in $C_\Th(H_1)$, $T\in\Cl_{H_1}(W)$. Take
$T^\ast=\alpha(\vp)_W(T).$ It is an $H_2$-closed congruence in
$\vp(W).$ We set: $\Phi(W/T)=\vp(W)/T^\ast.$ It is an object in
the category $C_\Th(H_2).$ Define further $\Phi$ on the morphisms.
Let the morphism $\sigma: W_1/T\to W_2/T_2$ be given in
$C_\Th(H_1).$ This $\sigma$ determines a commutative diagram
$$
\CD
 W_1   @> s >> W_2\\
@V\mu_{T_1} VV @VV\mu_{T_2} V\\
W_1/T_1 @>\sigma>> W_2/T_2
\endCD
$$
\noindent
 Here $s$ is not determined uniquely by $\sigma$, but it
induces $\sigma,$\ $\overline s=\sigma.$ $\mu_{T_1}$ and
$\mu_{T_2}$ are natural homomorphisms.

Let us consider the diagram
$$
\CD
 \vp(W_1 )  @> \vp(s) >> \vp(W_2)\\
@V\mu_{T_1^\ast} VV @VV\mu_{T_2^\ast} V\\
\vp(W_1)/T_1^\ast @>\Phi(\sigma)>> \vp(W_2)/T_2^\ast
\endCD
$$
\noindent
 where $\mu_{T_1^\ast}=\Phi(\mu_{T_1}),$\
$\mu_{T_2^\ast}=\Phi(\mu_{T_2)}).$ We want to define
$\Phi(\sigma)$ to make the diagram be commutative.

We want to check that $\vp(s)$ induces a morphism
$\vp(W_1)/T_1^*\to\vp(W_2)/T_2^*$.
Check first
that if  $w_1,w_2\in\vp(W_1)$ and $w_1T_1^\ast w_2,$ then $\vp(s)(w_1)T_2^\ast\vp(s)(w_2).$
Take $\rho^\ast=\rho_{\vp(W_1)}T_1^\ast.$ Let $\mu\rho^\ast\mu',$\
$\mu,\mu'\in\End(\vp(W_1)).$ Take further $\nu,\nu'\in\End(W_1)$
with $\vp(\nu)=\mu,$\ $\vp(\nu')=\mu'.$ As before, we have
$\nu\rho\nu',$ where $\rho=\rho_{W_1}(T_1).$ For every $w\in W_1,$
we have $w^\nu T_1w^{\nu'}.$ This means also that
$\mu_{T_1}\nu=\mu_{T_1}{\nu'}.$ Applying the initial diagram, we
get $w^{\nu s}T_2 w^{\nu's}.$ We use now that $\alpha$  and $\beta$
commute. The definition of geometrical
coordination of algebras (commutativity of $\alpha$ and $\beta$)
implies $w^{\vp(s\nu)}T_2^\ast w^{\vp(s\nu')}$ for every
$w\in\vp(W_2).$ We have:
$$w^{\mu\vp(s)}T_2^\ast w^{\mu'\vp(s)};$$
$$\vp(s)(w^\mu) T_2^\ast \vp(s)(w^{\mu'}).$$
We can find $w$ such  that $w^\mu=w_1,$\ $w^{\mu'}=w_2,$ which
leads to $\vp(s)(w_1)T_2^\ast\vp(s)(w_2).$ Hence, $\vp(s)$ induces
a homomorphism $\overline{\vp(s)}:\vp(W_1)/T_1^\ast\to
\vp(W_2)T_2^\ast.$ We set $\Phi(\sigma)=\Phi(\overline
s)=\overline{\vp(s)}.$

 We need also to check that this definition
of $\Phi(\sigma)$ does not depend on the choice of $s$ with
$\overline s=\sigma.$ Take $\mu_{T_2}s_1=\sigma\mu
_{T_1}=\mu_{T_2}s_2.$ For every $w\in W_1,$ we have
$\mu_{T_2}s_1(w)=\mu_{T_2}s_2(w)$ and $w^{s_1}T_2w^{s_2}.$ For
every $w\in\vp(W_1)$ we have $w^{\vp(s_1)}T_2^\ast w^{\vp(s_2)}.$
This follows from the commutativity with the function $\beta.$
Simultaneously, $\mu_{T_2^\ast}\vp(s_1)=\mu_{T_2^\ast}\vp(s_2).$
Take, further, an arbitrary
$$\overline w\in\vp(W_1)/T_1^\ast,\quad \overline w=w^{\mu_{T_2^\ast}},\quad w\in\vp(W_1).$$
Then $$\overline
w^{\overline{\vp(s_1)}}=\mu_{T_2^\ast}w^{\vp(s_1)}=\mu_{T_2^\ast}w^{\vp(s_2)}=\overline
w^{\overline{\vp(s_2)}}.$$ Hence, $\overline{\vp(s_1)}=\overline
{\vp(s_2)}=\Phi(\sigma).$ We have defined
$\Phi(\sigma):\Phi(W_1/T_1)\to\Phi(W_2/T_2)$ for an arbitrary
$\sigma: W_1/T_1\to W_2/T_2.$

Check that $\Phi$ carries   the multiplication of morphisms.

Let a commutative diagram
$$
\CD
   W_1    @>  s_1 >>  W_2  @>  s_2 >>  W_3\\
@V\mu_1=\mu_{T_1} VV @VV\mu_ 2=\mu_{T_2} V @VV \mu_ 3=\mu_{T_3} V\\
 W_1 /T_1  @>\sigma_1=\overline s_1>>  W_2/T_2   @>\sigma_2=\overline s_2>>  W_3/T_3
\endCD
$$
\noindent
\end{proof}
be given in $C_\Theta(H_1).$
Apply $\Phi:$
$$
\CD
  \varphi(W_1)    @>  \varphi(s_1) >>  \varphi(W_2) @>  \varphi(s_2) >>  \varphi(W_3)\\
@V\mu_1^\ast=\mu_{T_1^\ast}VV @VV\mu_ 2^\ast=\mu_{T_2^\ast} V @VV \mu_ 3^\ast=\mu_{T_3^\ast} V\\
\varphi(W_1) /T_1^\ast  @>\Phi(\ov {s_1})>>
 \varphi( W_2)/T_2^\ast   @>\Phi(\ov{s_2)}>>  \varphi(W_3)/T_3^\ast
\endCD
$$

{}From the first diagram we have $\ov{s_1s_2}=\ov{s_1}\ \ov{s_2};$
from the second one, we have
$$\ov{\varphi(s_1)\varphi(s_2)}=\Phi(\ov{s_1})\Phi(\ov{s_2})=\ov{\varphi(s_1s_2)}=\Phi(\ov{s_1s_2})=
\Phi(\ov{s_1}\ \ov{s_2}).$$ It is also clear that $\Phi(1)=1.$ Thus,
the functor $\Phi$ is built and it induces
$\varphi:\Theta^0\to\Theta^0.$ Similarly, we build
$\Psi:C_\Theta(H_2)\to C_\Theta(H_1)$ by $\psi,$ which also induces
$\psi:\Theta^0\to\Theta^0.$ It is left to check that $\Phi$ and
$\Psi$ give equivalence of categories.

We need to check that the product $\Psi\Phi=\Phi_0$ is an inner
autoequivalence of the category $C_\Theta(H_1),$ and
$\Phi\Psi=\Psi_0$ is an inner autoequivalence of the category
$C_\Theta(H_2).$

First fix $\psi\vp=\vp_0:\Theta^0\to\Theta^0$ and $\vp\psi=\psi_0:\Theta^0\to\Theta^0.$
These are inner autoequivalences.
Let $\vp_0$ relate to the isomorphism $s:1_{\Theta^0}\to\vp_0$ and $\psi_0$ is defined by the isomorphism
$s':1_{\Th^0}\to\psi_0.$
Extend these $s$ and $s'$ up to $S: 1_{C_\Theta(H_1)}\to \Phi_0$ and
$S':1_{C_\Theta(H_2)}\to\Psi_0.$
Since the autoequivalences $\vp_0$ and $\psi_0$ are inner, then for every $W\in\Ob\Theta^0$ and the congruence
$T$ in $W$ we have $\alpha(\vp_0)_W(T)=s_WT,$\ $\alpha(\psi_0)_W(T)=s_W'T.$
The isomorphism $s_W:W\to\vp_0(W)$ now induces the isomorphism
$$\bar s_W: W/T\to\Phi_0(W/T)=\vp_0(W)/T^\ast,$$
where $T^\ast=\alpha(\vp_0)_W(T).$ We have $\vp_0=\psi\vp.$
Further, use $\alpha(\psi\vp)_W(T)=\alpha
(\psi)_{\vp(W)}\alpha(\vp)_W(T).$ By  definition,
$$
\Phi(W/T)=\vp(W)/\alpha(\vp)_W(T);$$
$$\Phi_0(W/T)=\Psi\Phi(W/T)=\Psi(\vp(W)/\alpha(\vp)_W(T))=
\psi\varphi(W)/\alpha(\psi)_{\vp(W)}\alpha(\vp)_W(T)$$
$$=
\psi\vp(W)/\alpha(\psi\vp
)_W(T)=\vp_0(W)/\alpha(\vp_0)_W(T).
$$
Thus, $\Psi\Phi(W/T)=\psi\vp(W)/\alpha(\psi\vp)_W(T)$ and,
simultaneously, we have an isomorphism
$\ov{s_W}:W/T\to\Psi\Phi(W/T).$ Here $W/T$ is an arbitrary object
of the category $C_{\Th_1(H_1)}.$

Define now the function $S$ by the rule $S_{\overline W}: =
\ov{s_W}$, ${\overline W}=W/T$. Check that this defines the
isomorphism of functors $\ov S:
1_{C_\Theta(H_1)}\to\Psi\Phi=\Phi_0.$

Let the morphism $\sigma: W_1/T_1\to W_2/T_2$ be given in $C_\Theta(H_1),$ with the corresponding commutative
diagram
\noindent
$$
\CD
 W_1   @> \nu >> W_2\\
@V\mu_{T_1} VV @VV\mu{T_2} V\\
W_1/T_1 @>\sigma=\ov\nu>> W_2/T_2
\endCD
$$

We need to check that
$$
\CD
 W_1/T_1   @> \ov s_{\ov W_1} >>\Psi\Phi(W_1/T_1)\\
@V\sigma VV @VV\Psi\Phi(\sigma)V\\
W_2/T_2 @>\ov s_{\ov W_2} >> \Psi\Phi(W_2/T_2)
\endCD
$$
holds.
Here $\ov W_1=W_1/T_1,$\ $\ov W_2=W_2/T_2).$

Apply the functor $\Psi\Phi$ to the previous diagram.
$$
\CD
 \psi\vp(W_1 )  @> \psi\vp(\nu) >> \psi\vp(W_2)\\
@V\Psi\Phi(\mu_{T_1})VV @VV\Psi\Phi(\mu_{T_2}) V\\
\Psi\Phi(\ov W_1)  @>\Psi\Phi(\sigma)>> \Psi\Phi(\ov W_2)
\endCD
$$
Here $\Psi\Phi(\sigma)=\ov{\psi\vp(\nu)}.$
It does not depend on the choice of the representative $\nu.$
For $\nu: W_1\to W_2$ we have
$$
\CD
 W_1   @> s_{W_1} >>\psi\vp(W_1)\\
@V\nu VV @VV\psi\vp(\nu) V     \\
W_2 @>s_{W_2}>> \psi\vp(W_2)
\endCD
 $$
For $T_1$ from $W_1$ and $T_2$ from $W_2$ there hold the following
rules:
$$\ov{\nu}=\sigma: W_1/T_1\to W_2/T_2,\\$$
$$\ov{s_{W_1}}:W_1/T_1\to\Psi\Phi(W_1/T_1),\\$$
$$\ov{s_{W_2}}:W_2/T_2\to\Psi\Phi(W_2/T_2),\\$$
$$\Psi\Phi(\ov\nu)=\Psi\Phi(\sigma)=\ov{\psi\vp(\nu)}.$$
Now we check  commutativity of the  diagram
$$
\CD
 W_1/T_1   @> \ov{ s_{ W_1}} >>\psi\vp(W_1)/\alpha(\psi\vp)_{W_1}T_1)\\
@V\ov\nu VV @VV\ov{\psi\vp(\nu)}V\\
W_2/T_2 @>\ov {s_{ W_2}} >> \psi\vp(W_2/\alpha(\psi\vp)_{W_2}(T_2)
\endCD
$$
Rewrite the diagram in the following way:
$$
\CD
 \ov{W_1}   @> \ov{ s_{ W_1}} >>\psi\vp(W_1)/ \psi\vp(W_1)/s_{W_1}T_1\\
@V\ov\nu VV @VV\ov{\psi\vp(\nu)}V\\
 \ov{W_2}   @> \ov{ s_{ W_2}} >>\psi\vp(W_2)/ \psi\vp(W_2)/s_{W_2}T_2
\endCD
$$
But this diagram directly follows from the diagram for $\nu: W_1\to W_2$. 
We take $\ov w_1 \in \ov W_1 $ and act according to the rules
above. Thus, we have checked coordination of the function $S:
1_{C_\Theta(H_1)}\to\Psi\Phi$ with the morphisms of the category
and $S: 1_{C_\Theta(H_1)}\to\Psi\Phi$ is an isomorphism of
functors. We repeat the same  for $\vp\psi$ and $s':
1_{\Theta^0}\to\varphi\psi,$ thus coming to the isomorphism
$S':1_{C_\Theta(H_2)}\to\Phi\Psi.$ Correctness of the equivalence
of the categories $C_\Theta(H_1)$ and $C_\Theta(H_2)$ follows from
the fact that $\Phi$ induces $\vp,$\ $\Psi$ induces $\psi$  and,
by definition, $\Phi$ and $\Psi$ are coordinated with the natural
homomorphisms. The theorem is proved.

We will apply this theorem for the cases of the varieties
 $\Com$-$P$,  $\Ass$-$P$, and $Lie$-$P$.

\section{Automorphisms and Autoequivalences of Categories of Free Algebras of Varieties}\label{AutomAutoEq}

\subsection{The general problem and relation to main problems}

We are interested in automorphisms and autoequivalences of
categories of the form $\Theta^0$, where $\Theta^0$ is a variety
of algebras. The form of such automorphisms and autoequivalences
determines the peculiarities  of the similarity and coordination
relations. We have already seen that if all automorphisms of the
category $\Theta^0$ are inner, then all autoequivalences are inner
as well, and, thus for such $\Theta$ , the following conditions
are equivalent:
\begin{enumerate}
\item[1.] The algebras $H_1$ and $H_2$ in $\Theta$ are geometrically similar.
\item[2.] They are geometrically equivalent.
\item[3.] They are geometrically coordinated.
\end{enumerate}

Define further semi-inner automorphisms. We consider them in
general situation.

Let $\Theta$ be an arbitrary variety of algebras and $G$ an
algebra in $\Theta.$ Consider a new variety, denoted by
$\Theta^G.$ Define first the category $\Theta^G.$ Its objects have
the form $h: G\to H,$ where $H$ is an algebra in $\Theta$ and $h$
is a morphism in $\Theta.$ We call such objects {\it $G$-algebras}
in $\Theta,$ and denote them by $(H,h).$ The morphisms in
$\Theta^G$ are represented by commutative diagrams in $\Theta:$

$$
\CD
G @>h>> H\\
@. @/SE/ h'  // @VV \mu V\\
@. H' \\
\endCD
$$
An  algebra  $(H,h)$ is called {\it faithful} if $h$ is an
injection. We consider elements of the algebra $G$ as nullary
operations and add them to the signature of the variety $\Theta,$
thus gaining the variety $\Theta^G.$ For every set $X$ a free
algebra $W=W(X)$ in $\Theta^G$ is represented as free product
$$i_G:G\to G\ast W_0(X)=W(X),$$
where $W_0=W_0(X)$ is a free algebra in $\Theta$ over $X,$ ,
$i_G$ is an embedding related to free multiplication. Here $i_G$
turns out to be an injection.

In the category $\Theta^G$, along with its morphisms, consider
also semimorphisms. They are represented by  diagram
$$
\CD
 G   @> h>>H\\
@V\sigma VV @VV\nu V\\
G @>h' >> H'
\endCD
$$
where $\sigma$ is an endomoprhism of the algebra $G.$ We consider
a semimorphism as a pair $(\sigma,\nu),$ while a morphism is a
pair $(1,\nu).$ We consider semi-isomorphisms and
semi-automorphisms  for the objects from $\Theta^G.$

Let us pass to the category $(\Theta^G)^0$ of all free algebras $W=W(X)$ in $\Theta^G$
with finite $X.$

Define semi-inner automorphisms of this category.

\begin{defn}\label{defn7}
The automorphism $\vp:(\Theta^G)^0\to(\Theta^G)^0$ is called
semi-inner if there exists a semi-isomorphism of functors
$(\sigma,s):1_{(\Theta^G)^0}\to\vp$ with the automorphism $\sigma
$ of the algebra $G.$
\end{defn}

This means that for every object $W$ of the category
$(\Theta^G)^0$ a semi-isomorphism $(\sigma,s_W):W\to\vp(W)$ is
fixed, and for every morphism $\nu: W_1\to W_2$ we have
$$
\CD
 W_1   @> (\sigma,s_{W_1}) >> \vp(W_1)\\
@V\nu VV @VV \vp(\nu) V \\
W_2 @>(\sigma, s_{W_2})>>  \vp(W_2)
\endCD
$$
Here
 $\vp(\nu)=(\sigma,s_{W_2})(1,\nu)(\sigma^{-1},s_{W_1}^{-1})=(1,s_{W_2}\nu s_{W_1}^{-1})$
  is a morphism of
the category $(\Theta^G)^0.$

All semi-inner automorphisms of the category $(\Theta^G)^0$
constitute a subgroup in $\Aut(\Theta^G)^0,$ containing the
invariant subgroup $\Int(\Theta^G)^0.$

Varieties of algebras $\Ass$-$P$ and $\Com$-$P$ are varieties of
$\Theta^G$ type. Here $\Theta$ is the variety of associative
rings, with the unit in the first case, and $\Theta$ is the
variety of commutative and associative  rings with the unit in the
second case, where $G=P$ is a field. The first case assumes  that
embeddings $h: P\to H$ are embeddings into the center of the ring
$H.$

Consider the corresponding semimorphisms $(\sigma,s): H\to H'.$
Here $s: H\to H'$ is a homomorphism of rings and $s(\lambda
a)=\lambda^\sigma s(a),$ \ $\lambda\in P,$\ $a\in H$.

We consider semimorphisms also in the category of modules
$\Mod$-$K$ and the category of Lie algebras over a field. These
varieties are not varieties of the $\Theta^G$ type. However,
semi-inner automorphisms are naturally defined for the categories
$(\Mod$-$K)^0$ and $(\Lie$-$P)^0.$

Let us quote results from \cite{23}.
\begin{enumerate}
\item[1.] If $\Theta=\Grp$ is a variety of all groups, then all automorphisms of the category $\Theta^0$ are
inner.
\item[2.] If $\Theta$ is a variety of all semigroups, then the group
$\Aut(\Theta^0)$ is a direct product of
the group $\Int(\Theta^0)$ and a cyclic group of order two.
\item[3.] All automorphisms of the category $(\Com$-$P)^0$ are semi-inner.
\item[4.] If the ring $K$ is left-noetherian, then all automorphisms of
the category $(\Mod$-$K)^0$ are semi-inner.
\item[5.] If $F$ is a free group of finite rank, then all automorphisms of
the category $\Grp^F$ are semi-inner.
\end{enumerate}
This list    of results can be accomplished by the result on Lie algebras
(see Theorem 10).

Correspondingly, autoequivalences of categories are described in all these cases.

Recall that every autoequivalence $\vp$ of the category $\Theta^0$
has the form $\vp=\vp_0\zeta=\zeta\psi_0$, where $\zeta$ is an
automorphism and $\vp_0$ and $\psi_0$ are inner. See also [10],
[11].

\subsection{Semi-inner automorphisms and autoequivalences}

The definitions are already given above; now we consider some details.
We return to the situation $\Theta^G.$

To every automorphism $\sigma$ of the algebra $G$, we construct
the corresponding semi-inner automorphism $\hat\sigma$ of the
category $(\Theta^G)^0.$ For every $W=W(X)=G\ast W_0(X)$ we have
two embeddings
$$i_G\sigma: G\to G\ast W_0,$$
$$i_{W_0}: W_0\to G\ast W_0,$$
This gives the corresponding endomorphism in $\Theta$
$$\sigma_W: G\ast W_0\to G\ast W_0.$$
We have also an inverse endomorphism $\sigma^{-1}_W$ and, hence, $\sigma_W$ is an automorphism in $\Theta.$

It is easy to understand that the commutative diagram
$$
\CD
 G   @> i_G>>W\\
@V\sigma VV @VV\sigma_W V\\
G @>i_G >> W
\endCD
$$
takes place, and, thus, the pair $(\sigma,\sigma_W)$ defines the
semi-automorphism of the algebra $W.$

Let, further, $(\sigma,s)$ be an arbitrary pair, such that  for
every $W\in\Ob(\Theta^G)^0$ a semi-automorphism $(\sigma,s_W):W\to
W$ be fixed. The pair $(\sigma,s)$ defines a semi-inner
automorphism $\widehat{(\sigma,s)} $ of the category
$(\Theta^G)^0.$ It does not change objects and for every $\nu:
W_1\to W_2$ we have $\widehat{(\sigma,s)}(\nu)=s_{W_2}\nu
s_{W_1}^{-1}.$ In particular, if always $s_W=\sigma_W,$ then a
semi-inner automorphism of the category $(\Theta^G)^0,$ denoted by
$\hat\sigma,$ corresponds to the automorphism $\sigma\in\Aut(G).$

Consideration of the pair $(1,s)$ leads to the inner automorphism $\hat s$ of the category
$(\Theta^G)^0.$

Let us now treat semi-inner autoequivalences. They are defined
exactly in the same way as semi-inner automorphisms. An
autoequivalence $\vp:(\Theta^G)^0\to(\Theta^G)^0$ is semi-inner if
it is given by a semi-isomorphism of functors
$(\sigma,s):1_{(\Theta^G)^0}\to\vp.$

Show that every such $\vp$ can be represented as
$\vp=\hat\sigma\vp_1=\vp_0{\hat\sigma},$ where
$\vp_0$ and $\vp_1$ are inner autoequivalences.

Let a semi-inner autoequivalence $\vp$ be given by $(\sigma,s):
1_{\Theta^G)^0}\to\vp.$ For every $W\in\Ob(\Theta^G)^0$ we have a
semi-isomorphism $(\sigma, s_W): W\to \vp(W).$ Consider also
$(\sigma,\sigma_W): W\to W$ and $(\sigma^{-1},\sigma_W^{-1}): W\to
W.$

Take a product
$$(\sigma,s_W)(\sigma^{-1},\sigma_W^{-1})=(1,s_W\sigma_W^{-1}): W\to\vp(W).$$
Denote $s_W'=s_W\sigma_W^{-1}.$
We have an isomorphism $s_W':W\to \vp(W).$

Consider a function $s'$ defined by the rule $s_W'=s_W\sigma_W^{-1}.$
The function $s'$ determines the inner
autoequivalence $\vp_0: (\Theta^G)^0\to(\Theta^G)^0$ acting on
the objects as $\vp $ does: $\vp_0(W)=\vp(W).$
We have $(\sigma,s_W)=(1,s_W')(\sigma,\sigma_W)$ and, correspondingly, $\vp=\vp_0\hat \sigma.$
Similarly, we define the decomposition $\vp=\hat\sigma\vp_1.$
The same considerations are applicable to automorphisms.

\subsection{Application}

The following proposition was proved in \cite{7}.

\begin{prop}\label{prop13}
If the algebras $H_1$ and $H_2$ are geometrically similar and
their similarity is defined by the semi-inner automorphism, then
there exists an algebra $H$ which is semi-isomorphic to the
algebra $H_1$ and geometrically equivalent to the algebra $H_2.$
\end{prop}

The existence of such $H$ means that $H_1$ and $H_2$ are similar.

We want to prove also a similar proposition for the relation of
geometrical coordination of algebras, but first let us make an
auxiliary remark.

Given a $G$-algebra $(H,h)$ and $\sigma\in\Aut(G),$ build a new $G$-algebra $(H_1,h_1)$, \ $H_1=H,$
keeping in mind the commutative diagram
$$
\CD
 G   @> h_1>>H_1\\
@V\sigma  VV @VV\mu=1 V\\
G @>h >> H
\endCD
$$
Here $h_1=h\sigma$ and the algebras $(H,h)$ and $(H_1,h_1)$ are
semi-isomorphic.

\begin{prop}\label{prop14} (See [7]).
The algebras $(H,h)$ and $(H_1,h_1)$ are geometrically similar, and their similarity is defined by the
automorphism $\hat\sigma:(\Theta^G)^0\to(\Theta^G)^0.$
\end{prop}

Recall that we consider the situation $\Var(H)=\Theta^G.$ In this
case, it is easy to check that $\Var(H_1)=\Theta^G$ holds.
Besides, the algebras $(H,h)$ and $(H_1,h_1)$ are faithful
$G$-algebras.

\begin{prop}\label{prop15}
Let $\Var(H_1)=\Var(H_2)=\Theta^G$ and let the $G$-algebras $H_1$
and $H_2$ be geometrically coordinated by the semi-inner
autoequivalence $(\vp,\psi).$ Then there exists a $G$-algebra $H,$
semi-isomorphic to $H_1$ and geometrically equivalent to $H_2.$ In
particular, $H_1$ and $H_2$ are geometrically similar.
\end{prop}

\begin{proof}
We use  \propref{prop12}. Let $\vp$ and $\psi$ {\it be related to}
the automorphism $\sigma$ of the algebra $G.$ This means that we
can proceed from the decomposition $\vp=\vp_0\hat\sigma,$\
$\psi=\hat \sigma^{-1}\psi_0,$ where $(\vp_0,\psi_0)$ is an inner
autoequivalence of the category $\Theta^0$ and
 the
automorphism $\hat\sigma$ of the category $(\Theta^G)^0$
corresponds to  $\sigma.$
By the given
$\sigma\in\Aut(G)$ take the $G$-algebra $H,$ semi-isomorphic to
the algebra $H_1,$ such that $\hat\sigma$ defines similarity of
the algebras $H_1$ and $H$ (see \propref{prop14}). According to
\propref{prop12}, the pair $(\vp_0,\psi_0)$ defines coordination
of the algebras $H$ and $H_2.$ Since the pair $(\vp_0,\psi_0)$ is
inner, the algebras $H$ and $H_2$ are geometrically equivalent.

The proposition is proved.
\end{proof}

\section{Varieties $\Com$-$P,$ $\Ass$-$P$ and $\Lie-P$}\label{VarComPAssP}

\subsection{$\Theta= \Com$-$P$}

assume the field $P$ is infinite. The variety $\Com$-$P$ is
noetherian and is generated by each of its algebras. Two algebras
$H_1$ and $H_2$ are geometrically equivalent, if they have the
same quasi-identities. Besides, semimorhisms are naturally defined
in $\Com$-$P$. It is proved in \cite{7} that every automorphism of
the category $(\Com$-$P)^0$ is semi-inner. Then every
autoequivalence of this category is semi-inner as well. Now,
taking into account the previous considerations, we can formulate

\begin{thm}\label{thm8}
Let $H_1$ and $H_2$ be two algebras in $\Theta=\Com$-$P$.
Then the following conditions are equivalent.
\begin{enumerate}
\item[1.] The categories $K_\Theta(H_1)$ and $K_\Theta(H_2)$ are correctly isomorphic.
\item[2.] These categories are correctly equivalent.
\item[3.] There exists an algebra $H\in\Theta$ such that $H_1$ and $H$ are semi-isomorphic, and $H$ and $H_2$
have the same quasi-identities.
\end{enumerate}
\end{thm}

\subsection{$\Theta=\Ass$-$P$}

First of all, we are interested in automorphisms of the category
$\Theta^0=(\Ass$-$P)^0.$ In every category we can consider inner
automorphisms. This category $\Theta^0$ has also semi-inner
automorphisms. They are defined according to the general approach
from \secref{IsoEq}. Let us do it directly.

If $H_1$ and $H_2$ are two associative algebras over $P,$ then their
semimorphism $H_1\to H_2$ is given by the
pair $(\sigma,\nu),$ where $\sigma\in\Aut P,$ and $\nu: H_1\to H_2$ is a homomorphism of rings.
Here, if $\lambda\in P$ and $a\in H,$ then $\nu(\lambda a)=\lambda^\sigma\nu(a).$

In \secref{IsoEq} for every $G$-algebra $(H,h)$ and every
$\sigma\in\Aut(G)$, we considered the $G$-algebra $(H_1,h_1)$ with
$H$ and $H_1$ coinciding in $\Theta,$ and $h=h_1\sigma.$ The
$G$-algebras $H$ and $H_1$ are semi-isomorphic. Now we reproduce
this construction in $\Ass$-$P.$

Let $H$ be an associative algebra over the field $P.$
The embedding $h: P\to H$ is defined by the rule $\lambda^h=\lambda\cdot 1$
for every $\lambda\in P.$
Then $\lambda\cdot a=\lambda^h\cdot a$ for every $a\in H.$

Take a new algebra, denoted by $H^\sigma,$ for the given $\sigma\in\Aut(P).$
We set: $H$ and $H^\sigma$ coincide as rings, and we change multiplication by a scalar
$$\lambda\circ a=\lambda^{\sigma^{-1}} \cdot a=(\lambda^{\sigma^{-1}})^h\cdot a
=h\sigma^{-1}(\lambda)\cdot a=\lambda^{h_1}\cdot a.$$ Thus, the
algebra $H^\sigma$ is $(H_1,h_1)$ in the sense of general
construction. The identical transformation $H\to  H^\sigma$
determines semi-isomorphism of the algebras $H^\sigma$ and $H.$
Every semi-isomorphism can be decomposed in such a
semi-isomorphism and isomorphism.

We call the algebra $H^\sigma$ a $\sigma$-twisted algebra with respect to $H.$
It is checked that if $\Var(H)=\Theta,$ then $\Var(H^\sigma)=\Theta$ as well.

Consider further free algebras $W=W(X) $ in $\Theta=\Ass$-$P,$ cf
[9]. Denote by $S(X)$ a free monoid over $X$ and $S_0(X)$ a free
semigroup over $X.$ The algebra $W=W(X)$ is a semigroup algebra
$PS(X).$ Every element of $W(X)$ is uniquely represented in the
form
$$w=\lambda_0+\lambda_1u_1+\dots+\lambda_ku_k,\quad \lambda\in P,\ u\in S_0(X).$$
For every $\sigma\in\Aut(P)$ denote by $\sigma_W: W\to W$ a mapping, defined by the rule
$$\sigma_W(w)=\lambda_0^\sigma+\lambda_1^\sigma u_1+\dots +\lambda_k^\sigma u_k.$$
Here $\sigma_W$ is an automorphism of rings and the pair
$(\sigma,\sigma_W)$ defines semi-automorphism of the algebra $W.$

This definition corresponds to the general definition given above.

Denote by $\ov \sigma$ a function, choosing $\ov\sigma_W=\sigma_W$
for every $W.$ The pair $(\sigma,\ov\sigma)$ defines a semi-inner
automorphism $\hat\sigma$ of the category $\Theta^0.$ Here
$\hat\sigma$ does not change objects, and for every $\nu: W_1\to
W_2$ we have
$$\hat\sigma(\nu)=\sigma_{W_2}\cdot\nu\cdot\sigma_{W_1}^{-1}: W_1\to W_2.$$
An arbitrary semi-inner automorphism $\vp$ of the category
$\Theta^0$ is defined by the semi-isomorphism of the functors
$$(\sigma,s): 1_{\Theta^0}\to \vp.$$
Such a $\vp$ is represented as $\vp=\vp_0\hat\sigma=\hat\sigma\vp_1,$ where $\vp_0$ and $\vp_1$  are inner
automorphisms, and $\vp_0(W)=\vp(W)=\vp_1(W)$ for every $W.$

We have also semi-isomorphism
$$(\sigma,s_W):W\to\vp(W).$$
For $\nu:W_1\to W_2$ we have $\vp(\nu)=s_{W_2}\cdot\nu\cdot s_{W_1}^{-1}.$

Consider further a mirror automorphism of the category $(\Ass$-$P)^0.$
This notion relates to the idea of antimorphism in the category $\Ass$-$P.$
First consider antihomomorphisms  of semigroups.

A mapping of semigroups $\mu: S_1\to S_2$ is called an {\it
antihomomorphism } if $\mu(ab)=\mu(b)\mu(a),$\ $a,b\in S_1.$

Let now $S=S(X)$ be a free semigroup.
For every $u=x_{i_1}\dots x_{i_n}$ in $S$ take $\ov u=x_{i_n}\dots x_{i_1}.$
Then the transition $u\to\ov u$ is an antiautomorphism of the semigroup $S.$
Indeed, let $u=x_{i_1}\dots x_{i_n},$\ $v=x_{j_1}\dots x_{j_m}.$
Then
$$\ov{u\cdot v}=\ov{x_{i_1}\dots x_{i_n}x_{j_1}\dots x_{j_m}}=x_{j_m}\dots x_{j_1}x_{i_n}\dots x_{i_1}=\ov
v\cdot \ov u.$$
If now $H_1,$ $H_2$ are associative algebras over the field $P,$ then the mapping $\mu: H_1\to H_2$ is an
antihomomorphism of algebras if  $\mu$ is correlated with   addition and multiplication by a scalar, and
$\mu(ab)=\mu(b)\cdot \mu(a)$ for $a,b\in H_1.$

For an arbitrary algebra $H$ take an opposite algebra $H^\ast.$
The sets $H$ and $H^\ast$ coincide, $H$ and $H^\ast$ coincide also as vector spaces, but multiplication in
$H^\ast$ is defined by the rule $a\circ b=b\cdot a.$
An identical mapping $H\to H^\ast$ here is an antiisomorphism of algebras.

Let now $W=W(X)=PS(X)$ be a free associative algebra.

For every  its element $w=\lambda_0+\lambda_1u_1+\dots+\lambda_ku_k$ take $\ov w=\lambda_0+\lambda_1\ov u_1
+\dots +\lambda_k\ov u_k,$ and show that the transition $w\to \ov w$ is an antiautomorphism of the algebra $W.$

Given $w_1=\alpha_0+\alpha_1u+\dots+\alpha_ku_k$ and $w_2=\beta_0+\beta_1v_1+\dots+\beta_\ell v_\ell,$ we have
$$w_1w_2=\sum_{i,j}\alpha_i\beta_j u_i v_j,$$
$$\ov{w_1w_2}=\sum_{i,j}\alpha_i\beta_j\ov{u_iv_j}=\sum_{i,j}\beta_j\alpha_i\ov v_j\cdot\ov u_i=\ov w_2\cdot
\ov w_1.$$

Correlation with  addition and multiplication by a scalar are also evident.

Now we consider the mirror automorphism of the category $\Theta^0=(\Ass$-$P,)^0$ denoted by $\dl$.
This $\dl$ does not change objects.
Let the homomorphism $\nu:W_1=W(X)\to W(Y)=W_2$ be given.
Define $\dl(\nu):W_1\to W_2$ by $\dl(\nu)(x)=\ov{\nu(x)}$ for every $x\in X.$
Further we need additional calculations.

Let $u=x_{i_1}\dots x_{i_n}\in S_0(X).$
Consider
$$\dl(\nu)(u)=\dl(\nu)(x_{i_1})\dots
\dl(\nu)(x_{i_n})=\ov{\nu(x_i)}\dots\ov{\nu(x_{i_n})}=\ov{\nu(x_{i_n})\dots
\nu(x_{i_1})}=\ov{\nu(\ov{u})}.$$
If now $w=\lambda_0+\lambda_1u_1+\dots+\lambda_ku_k\in W(X),$ then
$$\dl(\nu)(w)=\lambda_0+\lambda_1\dl(\nu)(u_1)+\dots+\lambda_k\dl(\nu)(u_k)=\lambda_0+\lambda_1\ov{\nu(\ov{u}_1)
+\dots+\lambda_k\nu(\ov{u}_k)}\\$$
$$=\ov{\nu(\lambda_0+\lambda_1\ov{u}_1+\dots+\lambda_k\ov{u}_k)}=\ov{\nu(\ov{w})}.$$
Hence $\dl(\nu)(w)=\ov{\nu(\ov {w})}.$
Assume now that $\nu=\nu_1: W_1\to W_2$ and $\nu_2=W_2\to W_3$ are given.
Check that $\dl(\nu_2\nu_1)=\dl(\nu_2)\dl(\nu_1).$
Take an arbitrary $x\in X.$
Then
$$\dl(\nu_2)\cdot
\dl(\nu_1)(x)=\dl(\nu_2)(\overline{\nu_1(x)})=\overline{\nu_2(\overline{\overline{\nu_1(x)}})}=\ov{\nu_2\nu_1(x)}
= \dl(\nu_2\nu_1)(x).$$

It is also clear that $\dl(1)=1,$ and, thus, $\dl: \Theta^0\to\Theta^0$ is a functor.
Since $\dl^2=1_{\Theta^0},$ then $\dl$ is an automorphism.

Here $\dl$ is not inner and is not semi-inner, but is quasi-inner.
Besides, if $\widetilde{\Int}(\Theta^0)$ is a subgroup in
$\Aut(\Theta^0),$ consisting of semi-inner automorphisms, then
$\dl$ belongs to the normalizer of this subsemigroup.

Denote by $\eta$ a function, giving an antiautomorphism $\eta_W$ of the algebra $W$ by $\eta_W(w)=\ov w$ for
every $W=W(X).$
Show that $\dl(\nu)=\eta_{W_2}\cdot\nu\cdot\eta_{W_1}^{-1}$ holds for every $\nu: W_1\to W_2.$
Take an arbitrary $x\in X,$\ $W_1=W(X).$
Then
$$\eta_{W_2}\cdot\nu\cdot
\eta_{W_1}^{-1}(x)=\eta_{W_2}\cdot\nu(x)=\eta_{W_2}(\nu(x))=\ov{\nu(x)}=\dl(\nu)(x),$$
for every $x\in X.$ Hence,
$\dl(\nu)=\eta_{W_2}\cdot\nu\cdot\eta_{W_1}^{-1}.$ We checked that
$\dl$ is quasi-inner in this sense.

\begin{prop}\label{prop16}
The automorphism $\delta$ belongs to the normalizer of the subgroup in $\Aut(\Theta^0),$ consisting of
semiinner automorphisms.
\end{prop}

\begin{proof}
Let, first $\vp$ be an inner automorphism, defined by the
isomorphism of functors $s: 1_{\Theta^0}\to\vp.$ We have
$\dl^2=1_{\Theta^0},$\ $\dl^{-1}=\dl.$ Consider $\dl\vp\dl$ and
apply it to $\nu: W_1\to W_2.$ Then
$$\dl\vp\dl(\nu)=\dl(\vp(\dl(\nu)))=\dl(s_{W_2}\dl(\nu)s_{W_1}^{-1})=\\
$$
$$=\dl(s_{W_2})\dl^2(\nu)\dl(s_{W_1})^{-1}=\dl(s_{W_2})\nu\dl(s_{W_1})^{-1}.
$$
Thus, $\dl\vp\dl$ is an inner automorphism, defined by the isomorphism $\dl(s): 1_{\Theta^0}\to
\dl\vp\dl^{-1},$ where $\dl(s)$ is a function defined by $\dl(s)_W=\dl(s_W).$
In the case of semigroups we have $\dl(s_W)=s_W,$ where $\dl\vp\dl=\vp.$
In our situation this is not true, and $\dl\vp\dl\ne\vp.$
Indeed, if $s_W(x)=\lambda_0+\lambda_1u_1+\dots+\lambda_k u_k,$ where all $u_i$ depend on many variables, then
$\dl(s_W)(x)=\ov{s_W(x)}\ne s_W(x).$

Let further $\sigma\in\Aut(P).$
Consider the automorphism $\hat\sigma$ of the category $\Theta^0.$
Show that $\hat\sigma$ and $\dl$ commute.
Proceed once more from $\nu: W_1\to W_2,$ and check that $\hat\sigma\dl(\nu)=\dl\hat\sigma(\nu).$
Let $W_1=W(X),$\ $x\in X.$
Then
$$\hat\sigma(\dl(\nu))=\sigma_{W_2}\dl(\nu)\dl_{W_1}^{-1};$$
$$\hat\sigma(\dl(\nu))(x)=\sigma_{W_2}\dl(\nu)\sigma_{W_1}^{-1}(x)=\sigma_{W_2}\dl(\nu)(x)=\sigma_{W_2}\overline{(
\nu(x))}.$$ Let
$\nu(x)=\lambda_0+\lambda_1u_1+\dots+\lambda_ku_k.$ Then
$$\sigma_{W_2}\dl(\nu)(x)=\sigma_{W_2}(\lambda_0+\lambda_1\ov u_1+\dots+\lambda_k\ov
u_k)=\lambda_0^\sigma+\lambda_1^\sigma\ov u_1+\dots+\lambda_k^\sigma\ov u_k.$$
Here all $u_i$ are elements of $S_0(Y),$\ $W_2=W(Y).$
Now
$$\dl\hat\sigma(\nu)(x)=\dl(\sigma_{W_2}\nu s_{W_1}^{-1})(x)
=\ov{\sigma_{W_2}\nu\sigma_{W_1}^{-1}(x)}=
$$
$$\ov{\sigma_{W_2}(\lambda_0+\lambda_1u_1+\dots+\lambda_ku_k)}\\
=\ov{\lambda_0^\sigma+\lambda_1^\sigma u_1+\dots+\lambda_k^\sigma u_k}=\lambda_0^\sigma+\lambda_1^\sigma \ov
u_1+\dots+\lambda_k^\sigma\ov u_k.
$$
The proposition is proved.
\end{proof}

\begin{cor*}\label{cor*}
If $\vp$ belongs to a subgroup generated by semi-inner
automorphisms and the automorphism $\dl$ then $\vp$ is either a
semi-inner automorphism, or $\vp=\vp_0\dl,$ where $\vp_0$ is a
semi-inner automorphism.
\end{cor*}

\begin{prob}\label{prob6}
Whether it is true that the group $\Aut(\Ass$-$P)^0$ is generated by semi-inner and mirror automorphisms?
\end{prob}

\begin{prob}\label{prob7}
Let $F=F(X)$ be a free non-commutative Lie algebra. Whether it is
true that every automorphism of the semigroup $\End F$ is
semi-inner?
\end{prob}

The similar result for the category of free Lie algebras is proved.

Let us pass to the geometrical problems.
For every free algebra $W$ consider its antiautomorphism $\eta_W: W\to W$,\ $\eta_W(w)=\ov w.$
It is clear that if $T$ is  an ideal in $W,$ then its image $\eta_W(T)=T^\ast$ is also an ideal, and $w\in T^\ast$ if $\ov
w\in T.$
Check that $\alpha(\dl)_W(T)=T^\ast.$

Take $\rho=\rho_W(T)$ and $\rho^\ast=\rho_W(T^\ast).$ Verify that
$\rho^\ast=\dl(\rho).$ Let $\nu\rho\nu'$ hold. For every $w\in W$
and $w_1=\ov w\in W$ we have $\nu(w)-\nu'(w)\in T;$
$$\dl(\nu)(w)-\dl(\nu')(w)=\ov{\nu(\ov w)}-\ov{v'(\ov w)}=\ov{\nu(\ov w)-\nu'(\ov w)}\in T^\ast.$$
Therefore, $\dl(\nu)\rho^\ast\dl(\nu').$

Let now $\mu\rho^\ast\mu'.$
Take $\mu=\dl(\nu),$\ $\mu'=\dl(\nu').$
We have $\dl(\nu)(w)-\dl(\nu')(w)=\ov{\nu(\ov w)-\nu'(\ov w)}\in T^\ast$ for every $w\in  W,$ in which case $\nu(\ov
w)-\nu'(\ov w)\in T,\nu\rho\nu'.$
The equality $\dl(\rho)=\rho^\ast$ is verified.
Further,
$$
\alpha(\dl)(T)=\tau_W(\dl(\rho_W(T))=\tau_W(\dl(\rho))=\tau_W(\rho^\ast)=T^\ast.$$

\begin{prop}\label{prop17}
Let the algebras $H_1$ and $H_2$ be antiisomorphic.
Then they are geometrically similar,  and similarity is defined by the automorphism $\dl:\Theta^0\to\Theta^0.$
\end{prop}

\begin{proof}
Let $\mu: H_1\to H_2$ be an antiisomorphism.
Consider the commutative diagram

$$
\CD
W @> \eta_W>> W\\
@V\nu VV @VV \nu' V\\
H_1^I @>\ov\mu>> H_2^I
\endCD
$$
where $\ov\mu$ is an antiisomorphism defined by the
antiisomorphism $\mu,$ and  $\nu,$ $\nu'$ are one-to-one
corresponding homomorphisms of algebras. Prove now that if $T$ is
an $H_1$-closed ideal, then $T^\ast$ is $H_2$-closed, and vice
versa.


An injection $W/T\to H_1^I$ can be substituted by a homomorphism $\nu: W\to H_1^I$ with the kernel $T.$
It is easy to see that $T$ is $\Ker(\nu)$ if and only if $T^\ast$ is $\Ker(\nu').$
Hence the embedding $W/T\to H_1^I$ defines the embedding $W/T^\ast\to H_2^I,$ and vice versa.

It is left to check that $\dl$ and the function $\beta$ commute.
It is done in the same way as for $\dl(\rho)=\rho^\ast.$
\end{proof}

We call the algebras $H_1$ ad $H_2$ {\it almost geometrically
equivalent} if there exists a sequence $H_1,H,H',H_2$ such that
$H_1$ and $H$ are antiisomorphic or isomorphic, $H$ and $H'$ are
semi-isomorphic, and $H'$ and $H_2$ are geometrically equivalent.

We can now state that if $\Var(H_1)=\Var(H_2)=\Ass$-$P$ and
Problem 6
 about automorphisms of the category
$(\Ass$-$P)^0$ is solved positively, then the following conditions are equivalent:
\begin{enumerate}
\item[1.] Categories $K_\Theta(H_1)$ and $K_\Theta(H_2)$ are correctly isomorphic.
\item[2.]  They are correctly equivalent.
\item[3.] $H_1$ and $H_2$  are almost geometrically  equivalent.
\end{enumerate}
This is the main conjecture.
Let us discuss the statement  in more detail .
Let the algebras $H_1$ and $H_2$ be coordinated by an autoequivalence $(\vp,\psi)$ of the category $\Theta^0.$
Assume that $\vp=\vp_0\hat\sigma\dl,$\ $\psi=\dl^{-1}\hat\sigma^{-1}\psi_0,$ where $(\vp_0,\psi_0)$ is an inner
autoequivalence.
For $\dl$ take an algebra $H,$ opposite to $H_1.$
The algebras $H$ and $H_1$  are similar in respect to $\dl,$ and $H$ and $H_2$ are coordinated with respect to
$(\vp_0\hat\sigma,\hat\sigma^{-1}\psi_0)$ (\propref{prop12}).
Take an algebra $H'$ by $H,$ which is $\sigma$-twisted with respect to $H.$
The algebras $H$ and $H'$ are similar with respect to $\hat \sigma,$
\ $H'$ and $H_2$ are coordinated with respect to $(\vp_0,\psi_0).$
Since $(\vp_0,\psi_0)$ is an inner autoequivalence, $H'$ and $H_2$ are geometrically equivalent.


\subsection{Variety of Lie algebras Lie-$P$}

The following theorem is proved in \cite{24}.

\begin{thm}\label{thm10}
Every automorphism of the category of free Lie algebras is
semi-inner. Every autoequivalence of this category is semi-inner
as well.
\end{thm}

Consider an application of this theorem.

For every Lie algebra $H$ and every automorphism $\sigma$ of
the field $P$ consider a Lie algebra $H^\sigma,$ coinciding
with $H$ as a ring, while the multiplication by a scalar
is defined the new rule:

$$\lambda\circ a=\lambda^{\sigma^{-1}}\cdot a;\quad \lambda a=
\lambda^\sigma\circ a.$$ The identity mapping $H\to H^\sigma$ is a
semi-isomorphism of algebras.

The following theorem takes place:

\begin{thm}\label{thm11}
Let $\Var(H_1)=\Var(H_2)=\Lie$-$P$.
Then the following conditions are equivalent:
\begin{enumerate}
\item[1.] The categories $K_\Theta(H_1)$ and $K_\Theta(H_2)$ are isomorphic.
\item[2.] These categories are equivalent.
\item[3.] The algebra $H_1^\sigma$ is geometrically equivalent to the algebra $H_2$ for some $\sigma\in\Aut(P).$
\end{enumerate}
\end{thm}

\begin{proof}
Prove first that for any algebra $H$ the algebras $H$ and $H^\sigma$ are geometrically similar with respect to an
automorphism $\hat\sigma:\Theta^0\to\Theta^0;$\ $\Theta=\Lie$-$P.$

Define the automorphism $\hat\sigma.$ Let $W=W(X)$ be a free Lie
algebra over $P$ with finite $X$. Define for it a
semi-automorphism $\sigma_W:W\to W.$ Apply $\sigma_W$ to an
element $w\in W.$ We define the action of $\sigma_W$ inductively.
Set: $\sigma_W(x)=x$ for every $x\in X.$ If $w=w_1\cdot w_2,$ then
$\sigma_W(w)=\sigma_W(w_1)\cdot \sigma_W(w_2).$ Analogously, if
$w=w_1+w_2,$ then $\sigma_W(w)=\sigma_W(w_1)+\sigma_W(w_2).$ If,
finally, $w=\lambda w_1,$\ $\lambda\in P,$ then
$\sigma_W(w)=\lambda^\sigma\cdot \sigma_W(w).$ It can be verified
with the help of the suitable basis in $W$ that this definition is
correct. The pair $(\sigma,\sigma_W)$ determines a
semi-automorphism of the algebra $W.$

Set further: $\hat\sigma(W)=W$ for every $W\in\Ob\Theta^0$ and
$\hat\sigma(\nu)=\sigma_{W_2}\nu\sigma_{W_1}^{-1}$ for $\nu:
W_1\to W_2.$ This defines the semi-inner automorphism
$\hat\sigma:\Theta^0\to\Theta^0.$

We could not define here $\hat\sigma$ via the general approach, applied to varieties of $\Theta^G$ type, since the variety
$\Lie$-$P$ is not of such type.

Let us now link homomorphisms $W\to H$ and $W\to H^\sigma.$ Take
$\mu=\nu^\ast:W\to H^\sigma$ corresponds to $\nu: W\to H$ by the
rule $\mu(\sigma_W(w))=\nu(w),$\ $w\in W.$ Here
$\mu\sigma_W=\nu,$\ $\mu=\nu\sigma^{-1}_W,$\
$\mu(w)=\nu(\sigma_W^{-1}w).$ The mapping $\mu$ is coordinated
with the operations of the ring.

Check now that $\mu$ is a  homomorphism of algebras.
Indeed,
$$\mu(\lambda w)=\nu(\sigma_W^{-1}(\lambda w))=
\nu(\lambda^{\sigma^{-1}}\sigma_W^{-1}(w))\\
=\lambda^{\sigma^{-1}}\nu(\sigma_W^{-1}(w))=\lambda^{\sigma^{-1}}\mu(w)=\lambda\circ\mu(w).
$$

We have also: $w\in \Ker\nu$ if and only if $\sigma_W(w)\in\Ker\mu.$
If $A$ is a set of $H$-points, \ $A\subset\Hom(W,H),$ then a set $A^\ast$ of $H^\sigma$-points,\ $A^\ast\subset
\Hom(W,H^\sigma),$ corresponds to $A.$

Let now $T$ be an ideal in $W.$
Denote by $\sigma_WT$ an ideal in $W,$ consisting of all $\sigma_W(w),$\ $w\in T.$
It is clear now that
$$\capl_{\nu\in A}\Ker\nu=T\Leftrightarrow \capl_{\mu\in A^\ast}\Ker\mu=\sigma_WT.$$
This means that the ideal $T$ is $H$-closed if and only if the ideal $\sigma_WT$ is $H^\sigma$-closed.

Check that the transition $T\to T^\ast=\sigma_WT$ is coordinated with the function $\beta $ for $\varphi=\hat\sigma.$

Take algebras $W_1$ and $W_2$ in $\Ob\Theta^0.$
Let $T$ be an ideal in $W_2.$
Denote $\beta=\beta_{W_1,W_2}(T),$\ $\beta^\ast=\beta_{W_1,W_2}(T^\ast)$ for $T^\ast=\sigma_{W_2}(T).$
We need to check that
$$s\beta s'\Leftrightarrow \hat\sigma(s)\beta^*\hat\sigma(s')$$
for $s,s': W_1\to W_2.$
We have
$$\hat\sigma(s)=\sigma_{W_2}s\sigma_{W_1}^{-1},\quad \hat\sigma(s')=\sigma_{W_2}s'\sigma_{W_1}^{-1}.$$
Take an arbitrary $w\in W_1$ and consider a difference
$$\hat\sigma(s)(w)-\hat\sigma(s')(w)=\sigma_{W_2}(s(\sigma^{-1}_{W_1}(w))-s'(\sigma_{W_1}^{-1}(w)).$$
An arbitrary element in $w_1 \in W_1$ has the form
$w_1=\sigma_{W_1}^{-1}(w).$ Let now $s\beta s'$ take place. Then
$s(w_1)-s'(w_1)\in T.$ Hence $\hat \sigma(s)(w)-\hat(s')(w)\in
T^\ast,$ which gives $\hat\sigma(s)\beta^\ast\hat\sigma(s').$

It is also clear that if $\hat\sigma(s)(w)-\hat\sigma(s')(w)\in
T^\ast,$ then $s(w_1)-s'(w_1)\in T$ and, therefore, $s\beta s'$

Prove now that
$$\alpha(\hat\sigma)_W(T)=T^\ast=\sigma_WT.$$
It follows from considerations above that $\hat\sigma(\rho_W(T))=\rho_W(T^\ast).$
Applying $\tau_W,$ we get
$$\alpha(\hat\sigma)_W(T)=\tau_W(\hat\sigma(\rho_W(T))=\tau_W\rho_W(T^\ast)=T^\ast.$$
We have checked that there is a bijection
$$\alpha(\hat\sigma)_W: \Cl_H(W)\to\Cl_{H^\sigma}(W)$$
and the function $\alpha$ commutes with $\beta.$
This means that the automorphism $\hat\sigma$ determines similarity of algebras $H$ and $H^\sigma.$

Let us now finish the proof of the theorem.

The categories $K_\Theta(H_1)$ and $K_\Theta(H_2)$  are isomorphic if and only if $H_1$ and $H_2$ are similar.

The similarity of $H_1$ and $H_2$ is determined by some
automorphism $\varphi: \Theta^0\to \Theta^0.$ According to
\thmref{thm10}, an automorphism $\varphi$ is semi-inner and it can
be represented as $\varphi=\hat\sigma\varphi_0$ where $\varphi_0$
is an inner automorphism.

Let us pass to the algebra $H_1^\sigma.$
The algebras $H_1$ and $H_1^\sigma$ are similar in respect to $\hat\sigma.$
According to the similarity decomposition rule we conclude that $H_1^\sigma$ and $H_2$ are similar with respect to
$\varphi_0$ and, consequently, they are geometrically equivalent.
This leads to the equivalence of the first and the third conditions of \thmref{thm11}.
Equivalence of the second and the third connections is checked similarly.

\subsection{Acknowledgements}
The author is happy to thank B.Kunyavskii, R.Li\-pya\-nsky,
G.Mashevitzky, E.Plotkin, E.Rips, G.Zhitomiskii, N.Vavilov for the
stimulating discussions of the results. The final version of the
paper has been prepared in Jurmala, guesthouse "Allat". I am very
grateful to its director Galina Tuch and the stuff for the
wonderful hospitality and creative atmosphere.

\end{proof}

\end{document}
\end
\bye